\documentclass[preprint,12pt]{elsarticle}

\usepackage{graphicx}
\usepackage{amssymb}

\usepackage{lineno}
\usepackage{verbatim}
\usepackage{graphicx}
\usepackage{amsmath}
\usepackage{autobreak}
\usepackage{framed}
\usepackage{booktabs}
\usepackage{threeparttable}
\usepackage{geometry}
\biboptions{numbers,sort&compress}
\usepackage[section]{placeins}
\journal{Journal Name}

\begin{document}

\begin{frontmatter}


\title{Second-order computational homogenization of flexoelectric composites}



\author[1,2]{Xiaoying Zhuang\corref{mycorrespondingauthor}}
\author[1,2,3]{Bin Li \corref{mycorrespondingauthor}}
\author[2]{S.S. Nanthakumar}
\author[3]{Thomas Böhlke}

\address[1]{Department of Geotechnical Engineering, Tongji University, Siping Road 1239, 200092 Shanghai, China}
\address[2]{Institute of Photonics, Department of Mathematics and Physics, Leibniz University Hannover, 30167 Hannover, Germany}
\address[3]{Institute of Engineering Mechanics, Karlsruhe Institute of Technology (KIT), 76131 Karlsruhe, Germany}

\begin{abstract}
Flexoelectricity shows promising applications for self-powered devices with its increased power density. 
This paper presents a second-order computational homogenization strategy for flexoelectric composite. The macro-micro scale transition, Hill-Mandel energy condition, periodic boundary conditions, and macroscopic constitutive tangents for the two-scale electromechanical coupling are investigated and considered in the homogenization formulation.The macrostructure and microstructure are discretized using $C^1$ triangular finite elements. The second-order multiscale solution scheme is implemented using  ABAQUS with user subroutines. Finally, we present numerical examples including parametric analysis of a square plate with holes and the design of 
piezoelectric materials made of non-piezoelectric materials to demonstrate the numerical implementation and the size-dependent effects of flexoelectricity.
\end{abstract}




\end{frontmatter}

\section{Introduction}

Flexoelectricity-an electromechanical coupling between polarization and strain gradient-exists in all dielectric materials. Moreover, flexoelectricity is a size-dependent effect that becomes more significant in nanoscale systems \cite{majdoub2008enhanced}. These features make it different from piezoelectricity, which is induced polarization due to homogeneous strain. For an overview of flexoelectricity, readers are referred to the articles \cite{nguyen2013nanoscale, zubko2013flexoelectric, yudin2013fundamentals, wang2019flexoelectricity, arias2022emancipation}.

The flexoelectric effect at the microscale has been exploited to design flexoelectric composites (piezoelectric metamaterials \cite{mocci2021geometrically}, multiferroic composites \cite{zhang2016enhancing}, etc).  The porous micro-structured materials were designed to generate an ultrahigh flexoelectric effect \cite{zhang2021ultrahigh}.  Moreover, the flexoelectric effect has long been regarded as a desirable property for advanced nano-/micro-electromechanical systems (N/MEMS) due to its universality and excellent scaling effect \cite{wang2019flexoelectricity}. The homogenized electromechanical properties and multiscale modeling are very important, however not sufficient, for the analysis and design of the heterogeneous flexoelectric microstructure. A micromorphic approach for modeling the scale-dependent effects of flexoelectricity was presented by McBride et al \cite{mcbride2020modelling}. The asymptotic expansion method was applied to derive effective coefficients of flexoelectric rods \cite{guinovart2019asymptotic} and flexoelectric membranes \cite{mohammadi2014theory}. The classic Eshelby's formalism was used for the homogenization of piezoelectric nanocomposites without using piezoelectric materials \cite{sharma2007possibility}.

The first-order computational homogenization method is a powerful approach to assess micro-macro structure-property relations \cite{geers2016multiscale}, where only the first gradient of the macroscopic displacement field is included. However, the uniformity assumption and neglect of geometrical size effects are the two major disadvantages of the first-order method, which significantly limit its applicability \cite{geers2010multi}. To overcome the concerns of the first-order homogenization method, the second-order computational homogenization method \cite{kaczmarczyk2008scale, kouznetsova2004multi, lesivcar2014second, lesivcar2016large} has been developed, where the gradient of the macroscopic deformation gradient tensor is incorporated in the $C^{1}-C^{0}$ macro-micro scale transition. In this approach, $C^1$ continuous interpolations are applied at the macrolevel, while standard $C^0$ continuity is used for the discretization of the RVE at the microlevel, so an additional integral condition on the fluctuation field at microlevel is required \cite{lesivcar2014second}. 

The extension to second-order computational homogenization with the $C^{1}-C^{1}$ macro-micro scale transition has been proposed \cite{lesivcar2017two,  lesivcar2015multiscale}, where the higher-order stress and strain tensors exist at both scales, and the homogenized stiffness tensors are calculated with the static condensation procedure. Due to the influence of strain gradient on material exhibiting flexoelectricity, the second-order homogenization method of $C^{1}-C^{1}$ macro-micro scale transition should be conducted. However, the static condensation procedure is not suitable for the electromechanical coupling problems, so perturbation analysis is instead used in this study. The global $C^{1}$ continuity for the discretization is required, and options include isogeometric analysis (IGA) \cite{ghasemi2017level, nanthakumar2017topology}, mixed FEM \cite{deng2017mixed, mao2016mixed, poya2019family}, meshless methods \cite{he2019characterizing} and $C^1$  triangular elements \cite{yvonnet2017numerical, dasgupta1990higher}.

In this paper, we develop a second-order homogenization method for flexoelectric composites accounting for size-dependent behaviour. The homogenized electromechanical properties are computed by adopting the perturbation analysis. The $C^{1}-C^{1}$ micro-macro scale transition, periodic boundary conditions, and the Hill-Mandel energy condition are discussed. $C^1$  triangular element with eighteen degrees of freedom per node is chosen for the discretization because of the ease of imposing periodic boundary conditions. The second-order homogenization scheme is implemented with the secondary development of ABAQUS.

The remainder of the paper is organized as follows: Section 2 introduces the flexoelectricity theory; Section 3 presents the formulation of second-order computational homogenization; Section 4 describes the numerical implementation; Several numerical examples are illustrated in Section 5 and the conclusion derived from this work is presented in Section 6.

\section{Flexoelectricity theory}

Within the assumption of the linearized theory for dielectric solid, the electrical enthalpy density $H$ can be written as \cite{maranganti2006electromechanical, sharma2010piezoelectric, abdollahi2015revisiting, shen2010theory}. 
\begin{gather}
  H=\frac{1}{2}c_{ijkl}\varepsilon_{i j} \varepsilon_{k l}+\frac{1}{2}b_{ijklmn}g_{i jk} g_{lmn}-e_{ijk}E_{i}\varepsilon_{jk}-f_{ijkl}E_{i}g_{jkl}-\frac{1}{2} \kappa _{ij}E_{i}E_{j}
  \label{H}
\end{gather}

where $\varepsilon_{ij}$, $g_{jkl}$, $E_{i}$ are the strain, strain-gradient, and electric field, respectively. They are defined as $\varepsilon_{i j}=\frac{1}{2}(u_{i,j}+u_{j,i})$, $g_{jkl}=u_{j,kl}$, and $E_{i}=-\phi_{,i}$, where $\boldsymbol{u}$ and ${\phi}$ are displacement vector and electric potential, respectively. $c_{ijkl}$, $b_{ijklmn}$, $e_{ijk}$, $f_{ijkl}$ and $\kappa_{ij}$ are the elastic, strain-gradient elastic,  piezoelectric, flexoelectric and dielectric tensors, respectively. 

The electric polarization $\boldsymbol{P}$ is introduced via the constitutive equation as \cite{yudin2013fundamentals}
\begin{gather}
  P_i =\chi_{ij}E_{j}+e_{ijk}\varepsilon_{jk}+\mu _{ijkl}\frac{\partial \varepsilon_{jk}}{\partial  x_l} =\chi_{ij}E_{j}+e_{ijk}\varepsilon_{jk}+f_{ijkl}\frac{\partial^2 u_{j}}{\partial  x_k\partial  x_l}
  \label{Pi1}
\end{gather}
where $\boldsymbol{\chi}=\epsilon_{0}\boldsymbol{\chi}_{e}$,  $\boldsymbol{\chi}$, $\boldsymbol{\chi}_{e}$ and  $\epsilon_{0} $ are the clamped dielectric susceptibility, electric susceptibility tensor and the permittivity of a vacuum, respectively. $\boldsymbol{f}$ and $\boldsymbol{\mu}$ are the flexoelectric coefficients corresponding to different descriptions of strain-gradient, and these two tensors are indirectly related due to the symmetry of the strain tensor.

There is an alternative formulation of flexoelectricity, assuming the internal energy $U$ as a function of strain, strain gradient, and polarization \cite{mao2016mixed},
\begin{gather}
  U=\frac{1}{2}c_{ijkl}\varepsilon_{i j} \varepsilon_{k l}+\frac{1}{2}b_{ijklmn}g_{i jk} g_{lmn}+d_{ijk}P_{i}\varepsilon_{jk}+h_{ijkl}P_{i}\varepsilon_{jk,l}+\frac{1}{2} \alpha _{ij}P_{i}P_{j}
  \label{U}
\end{gather}
where $U=H+E_{i}D_{i}$, $d_{ijk}$, $h_{ijkl}$ and $\alpha_{ij}$ are the piezoelectric, flexoelectric, and reciprocal dielectric susceptibility tensors, respectively. It is worth noting that $\boldsymbol{e}$ and $\boldsymbol{\mu}$ are the coupling between the polarization and the strain and strain gradient respectively, while  $\boldsymbol{d}$ and $\boldsymbol{h}$ links the electric field to the strain and strain gradient respectively \cite{maranganti2006electromechanical}. The tensor $\boldsymbol{e}$, $\boldsymbol{\mu}$ and $\boldsymbol{\kappa }$ in Eqs. (\ref{H})-(\ref{Pi1})  are related to the tensor $\boldsymbol{d}$, $\boldsymbol{h}$ and $\boldsymbol{\alpha}$ of Eq. (\ref{U}) as \cite{sharma2010piezoelectric,mao2016mixed}
\begin{gather}
d_{ijk}=\alpha_{im}e_{mjk},  \quad  h_{ijkl}=\alpha_{im}\mu_{mjkl},  \quad \kappa_{ij}=\chi_{ij}+ \epsilon_{0}\delta _{ij}=1/\alpha_{ij}+\epsilon_{0}\delta _{ij}
\end{gather}

For isotropic materials, the electrical enthalpy density of Eq. (\ref{H}) can be written as \cite{mao2016mixed}

\begin{gather}
  H=\frac{1}{2}  \lambda \varepsilon_{ii} \varepsilon_{jj}+G\varepsilon_{ij} \varepsilon_{ij}+  \frac{1}{2}l^2[\lambda g_{jij} g_{kik} + G ( g_{ijk} g_{ijk} + g_{ijk} g_{jik})]  -e_{ijk}E_{i}\varepsilon_{jk}  \\ \nonumber
  - (f_{1}\delta_{ij}\delta_{kl} +f_{2} (\delta_{ik}\delta_{jl}+\delta_{il}\delta_{jk}))E_{i}g_{jkl}    -\frac{1}{2} \kappa _{ij}E_{i} E_{j}
\end{gather}

where $\lambda$ and $G$ are the Lamé parameters, $l$ is the intrinsic length of material, $f_{1}$ and $f_{2}$ are the two independent flexoelectric constants of flexoelectric coefficient $f_{ijkl}$.

The corresponding constitutive equations are  
\begin{gather}
  \sigma _{ij}=\frac{\partial H}{\partial \varepsilon_{i j}} =2G  \varepsilon_{ij}+\lambda \varepsilon_{kk}\delta _{ij}-e_{kij}E_{k}\\
   \tau _{ijk}=\frac{1}{2}(\frac{\partial H}{\partial g_{ijk}}+ \frac{\partial H}{\partial g_{ikj}})=\frac{1}{2}l^2 [ \lambda (g_{njn} \delta _{ik}+g_{nkn} \delta _{ij}) + G ( 2g_{ijk}+ g_{jik} +g_{kij})]+   \\ \nonumber
  f_{1} \delta _{jk}(-E_{i})+f_{2}(\delta _{ij}(-E_{k})+\delta _{ik}(-E_{j}))    \\ 
  D_{i}=\frac{\partial H}{\partial (-E _{i})}=\kappa _{ij}E_{j}+e_{ijk}\varepsilon_{jk}+   f_{1}g_{ijj}+f_{2}(g_{jij}+g_{jji})  =P_{i}+\epsilon _{0}E_{i} 
\end{gather}
where $\sigma _{ij}$, $\tau _{ijk}$ and $D_{i}$ are the stress, higher-order stress, and electric displacement, respectively. 

The governing equations are given as \cite{hu2010variational,mao2016mixed}
\begin{gather}
  (\sigma_{ij}-\tau _{ijk,k})_{,j}+f_{i}^{b}=0 \quad \text{in} \quad \Omega\\
    D_{i,i}=\rho ^{f}  \quad \text{in} \quad \Omega
\end{gather}
where electrostatic stress and higher-order electric displacement are neglected, and $f_{i} ^{b}$ and $\rho ^{f}$ are body force and free charge on the solid body.\\

The boundary conditions are obtained as \cite{hu2010variational,deng2017mixed,mao2016mixed} 
\begin{gather}
  u_{i}=\bar{u}_{i} \quad \text{on}  \quad \Gamma_{u}\\
  \phi=\bar{\phi} \quad \text{on} \quad \Gamma_{\phi}\\
  (\sigma_{ij}-\tau _{ijk,k})n_{j}+(D_{l}n_{l})\tau _{ijk}n_{k}n_{j}-D_{j}(\tau _{ijk}n_{k})=\bar{t}_{i} \quad \text{on}  \quad \Gamma_{t} \\
  \tau _{ijk}n_{j}n_{k}=\bar{r}_{i} \quad \text{on}  \quad \Gamma_{r}\\
  u_{i,k}n_{k}=\bar{v}_{i} \quad  \text{on}  \quad \Gamma_{v}\\
 D_{i}n_{i}=-\omega \quad  \text{on}  \quad \Gamma_{D}
\end{gather}
where $D_{j}=(\delta_{jk}-n_{j}n_{k})\partial_{k}$ is the tangent gradient operator, and $\bar{t}_{i} $, $\bar{r}_{i}$, $\bar{v}_{i}$ and $\omega$ are surface traction, double traction, normal derivative of displacement, and surface charge density, respectively.

\section{Second-order computational homogenization}

\subsection{Macro-micro scale transition}

The displacement and electric potential field are represented by a Taylor series expansion \cite{kaczmarczyk2008scale}.
\begin{gather}
   \boldsymbol{u}_{m}=\boldsymbol{\varepsilon}_{M} \cdot \boldsymbol{x}+\frac{1}{2} \left [ \boldsymbol{g}_{M}  :  (\boldsymbol{x} \otimes \boldsymbol{x}) \right ]+\boldsymbol{r}_{u}  \\
   \phi _{m}=(-\boldsymbol{E}_{M}) \cdot \boldsymbol{x}+r_{\phi }
\end{gather} 
where $\boldsymbol{g}_{M} =(\boldsymbol{u} _{M}  \otimes \nabla \otimes \nabla)$, $\boldsymbol{x}$, $\boldsymbol{r}_u$ and $r_{\phi }$ are the microlevel spatial coordinate, the displacement microfluctuation field, and electric potential microfluctuation field, and the gradient of the electric field is ignored. The volume average of the strain and electric field at the microlevel are equal to the corresponding quantity at the macrostructural material point, as follows.
\begin{gather}
  \frac{1}{V}\int_{v} \boldsymbol{\varepsilon}_{m}dv=\frac{1}{V}\int_{v}(\boldsymbol{u} _{m} \otimes   \nabla )dv=
\boldsymbol{\varepsilon}_{M} + \frac{1}{V}\int_{v}(\boldsymbol{g}_{M}\cdot \boldsymbol{x} )dv + \frac{1}{V}\int_{V}(\boldsymbol{r}_{u}   \otimes  \nabla )dv \label{mstrain} \\
\frac{1}{V}\int_{v} \boldsymbol{g}_{m}dv =\frac{1}{V}\int_{v}(  \boldsymbol{u}_{m}  \otimes \nabla   \otimes \nabla )dv= \boldsymbol{g}_{M}+ 
\frac{1}{V}\int_{v}(    (  \boldsymbol{r}_{u}  \otimes  \nabla)  \otimes  \nabla  )dv \label{mgstrain}\\
\frac{1}{V}\int_{v} (-\boldsymbol{E}_{m})dv =\frac{1}{V}\int_{v}(\phi_{m}  \nabla)dv=
(-\boldsymbol{E}_{M}) + \frac{1}{V}\int_{v}( r_{\phi }  \nabla )dv \label{med}
\end{gather}
where the minor symmetry $(\boldsymbol{g}_{M})_{ijk}=(\boldsymbol{g}_{M})_{ikj}$ in Eq. (\ref{mstrain}) has been used. The second term on the right side of Eq. (\ref{mstrain}) $\frac{1}{V}\int_{v}(\boldsymbol{g}_{M}\cdot \boldsymbol{x} )dv$ should be zero, which indicates that the centroid of the RVE is at the origin of the coordinates. The following equations can also be obtained.
\begin{gather}
  \frac{1}{V}\int_{v}( \boldsymbol{r}_{u} \otimes \nabla )dv=\frac{1}{V}\int_{\Gamma }( \boldsymbol{r}_{u} \otimes \boldsymbol{n})d\Gamma =\boldsymbol{0}\\
  \frac{1}{V}\int_{v}(( \boldsymbol{r}_{u} \otimes \nabla) \otimes  \nabla))dv=\frac{1}{V}\int_{\Gamma }( ( \boldsymbol{r}_{u} \otimes \nabla)  \otimes \boldsymbol{n} )d\Gamma =\boldsymbol{0}\\
  \frac{1}{V}\int_{v}( r_{\phi } \nabla)dv=\frac{1}{V}\int_{\Gamma }( r_{\phi } \boldsymbol{n})d\Gamma =\boldsymbol{0}
\end{gather}

Therefore, the following boundary conditions can be obtained:
\begin{gather}
  \boldsymbol{r}_{uR}=\boldsymbol{r}_{uL},\quad
  \boldsymbol{r}_{uT}=\boldsymbol{r}_{uB}\\
  \boldsymbol{r}_{uR}\otimes \nabla=\boldsymbol{r}_{uL}\otimes \nabla,\quad
  \boldsymbol{r}_{uT}\otimes \nabla=\boldsymbol{r}_{uB}\otimes \nabla\\
  {r}_{\phi R}={r}_{\phi L},\quad
  {r}_{\phi T}={r}_{\phi B}
\end{gather}
where $R$, $L$, $T$ and $B$ represent right, left, top and bottom boundary of rectangular RVE, respectively.

\subsection{The Hill-Mandel energy condition}

Without considering the electric field gradient, the Hill-Mandel energy condition \cite{hill1963elastic} can be expressed as 
\begin{gather}
  \frac{1}{V} \int _{v} (\boldsymbol{\sigma}_m :\delta \boldsymbol{\varepsilon}_m+\boldsymbol{\tau}_m \vdots  \delta \boldsymbol{g}_m -\boldsymbol{D}_m \cdot \delta \boldsymbol{E}_m) dv=\boldsymbol{\sigma}_M :\delta \boldsymbol{\varepsilon}_M
  +\boldsymbol{\tau}_M \vdots \delta \boldsymbol{g}_M  -\boldsymbol{D}_M \cdot \delta \boldsymbol{E}_M
\end{gather}

This condition can be subdivided into a mechanical and an electrical part as \cite{schroder2012two}
\begin{gather}
  \frac{1}{V}\int _{v} (\boldsymbol{\sigma}_m :\delta \boldsymbol{\varepsilon}_m+\boldsymbol{\tau}_m \vdots \delta \boldsymbol{g}_m ) dv=\boldsymbol{\sigma}_M :\delta \boldsymbol{\varepsilon}_M +\boldsymbol{\tau}_M \vdots\delta \boldsymbol{g}_M \label{uhm}\\
  \frac{1}{V} \int _{v} (  -\boldsymbol{D}_m \cdot \delta \boldsymbol{E}_m) dv= -\boldsymbol{D}_M \cdot \delta \boldsymbol{E}_M \label{phihm}
\end{gather}

Substituting Eqs. (\ref{mstrain})-(\ref{mgstrain}) into the left side of Eq. (\ref{uhm}), we get
\begin{gather}
  \frac{1}{V} \int _{v} (\boldsymbol{\sigma}_m :(\delta \boldsymbol{\varepsilon}_M+ \delta \boldsymbol{g}_M \cdot \boldsymbol{x} + \delta \boldsymbol{r}_{u} \otimes \nabla )+\boldsymbol{\tau}_m \vdots (\delta \boldsymbol{g}_M +  ( \delta \boldsymbol{r}_{u} \otimes \nabla) \otimes  \nabla) ) dv\\ \nonumber
 = (\frac{1}{V}\int _{v} \boldsymbol{\sigma}_m dv):\delta \boldsymbol{\varepsilon}_M+ \frac{1}{V}\int _{v} (\boldsymbol{\sigma}_m :(\delta \boldsymbol{r}_{u} \otimes \nabla ))dv+ \\ \nonumber
 (\frac{1}{V}\int _{v} (\boldsymbol{\tau}_m +\boldsymbol{\sigma}_m \otimes \boldsymbol{x} )dv)\vdots \delta \boldsymbol{g}_M+\frac{1}{V} \int _{v}( \boldsymbol{\tau}_m \vdots ((\delta \boldsymbol{r}_{u} \otimes \nabla) \otimes  \nabla))dv
\end{gather}

To meet the Hill-Mandel condition (\ref{uhm}), the integral terms containing microfluctuations should vanish:
\begin{gather}
  \frac{1}{V}\int _{v} (\boldsymbol{\sigma}_m :(\delta \boldsymbol{r}_{u} \otimes \nabla ))dv+\frac{1}{V} \int _{v}( \boldsymbol{\tau}_m \vdots (( \delta \boldsymbol{r}_{u} \otimes \nabla) \otimes  \nabla))dv=0 \label{rsum}
\end{gather}

The second term of Eq. (\ref{rsum}) can be written as
\begin{gather}
  \frac{1}{V} \int _{v}( \boldsymbol{\tau}_m \vdots ((\delta \boldsymbol{r}_{u} \otimes \nabla) \otimes  \nabla))dv=  
  \frac{1}{V} \int _{v}(  (\boldsymbol{\tau}_m : (  \delta \boldsymbol{r}_{u} \otimes \nabla))  \cdot \nabla)dv-
  \frac{1}{V} \int _{v}( (\boldsymbol{\tau}_m   \cdot \nabla) : (\delta \boldsymbol{r}_{u} \otimes  \nabla) )dv \label{rsum2}
\end{gather}

Therefore, substituting Eq. (\ref{rsum2}) into Eq. (\ref{rsum}), we obtain:
\begin{gather}
  \frac{1}{V}\int _{v} ( (\boldsymbol{\sigma}_m- \boldsymbol{\tau}_m  \cdot \nabla) :(\delta \boldsymbol{r}_{u} \otimes  \nabla))dv+\frac{1}{V} \int _{v}(  (\boldsymbol{\tau}_m : (  \delta \boldsymbol{r}_{u} \otimes \nabla))  \cdot \nabla)dv=0  \label{rtsum1}
\end{gather}

Using Gauss theorem, Eq. (\ref{rtsum1}) can be transformed into the following form containing the integral over the boundary of RVE,
\begin{gather}
  \frac{1}{V}\int _{\varGamma } ( (\boldsymbol{\sigma}_m- \boldsymbol{\tau}_m  \cdot \nabla) \cdot \delta \boldsymbol{r}_{u}  \cdot  \boldsymbol{n})d\varGamma+ \frac{1}{V} \int _{\varGamma}(  (\boldsymbol{\tau}_m  : (\delta \boldsymbol{r}_{u} \otimes \nabla)) \cdot  \boldsymbol{n})d\varGamma=0  \label{rtsum2} 
\end{gather}

 where Eq. (\ref{rtsum2}) can be approximately satisfied by applying the periodic conditions which will be discussed in Section 4.2 due to the periodicity assumption and the suppression of the corner microfluctuations \cite{lesivcar2017two}. 

Similarly, the electrical part for the Hill-Mandel energy condition can be derived as follows. Substituting Eq. (\ref{med}) into the left side of Eq. (\ref{phihm}), we get
\begin{gather}
  \frac{1}{V} \int _{v} ( \boldsymbol{D}_m \cdot (\delta (-\boldsymbol{E}_M)+ (\delta r_\phi \nabla)) dv=
  ( \frac{1}{V} \int _{v} \boldsymbol{D}_m dv )\cdot \delta (-\boldsymbol{E}_M)+\frac{1}{V} \int _{v} ( \boldsymbol{D}_m \cdot (\delta r_\phi \nabla)) dv \label{phihm2}
  \end{gather}

The second term on the rigth-hand side of Eq. (\ref{phihm2}) should be zero:
\begin{gather}
  \frac{1}{V} \int _{v} ( \boldsymbol{D}_m \cdot (\delta r_\phi \nabla)) dv=\frac{1}{V} \int _{v}( (\delta r_\phi \boldsymbol{D}_m)  \cdot \nabla )dv
  -\frac{1}{V} \int _{v} ( \delta r_\phi ( \boldsymbol{D}_m  \cdot \nabla ) )=\frac{1}{V} \int _{\Gamma}( (\delta r_\phi \boldsymbol{D}_m)  \cdot \boldsymbol{n})d\Gamma=0 
\end{gather}
where $( \boldsymbol{D}_m  \cdot \nabla )=0$, as the free charge on the solid body is not considered in microscale. Eq. (37) can also be satisfied due to the periodicity assumption.

According to Eq. (31) and Eq. (36), the homogenized stresses and electric displacement can be derived as
\begin{gather}
\boldsymbol{\sigma }_M=\frac{1}{V}\int _{v}  \boldsymbol{\sigma }_m dv \\
\boldsymbol{\tau  }_M=\frac{1}{V}\int _{v}  (\boldsymbol{\tau  }_m+\boldsymbol{\sigma }_m \otimes  \boldsymbol{x }) dv\\
\boldsymbol{D}_M=\frac{1}{V} \int _{v} \boldsymbol{D}_m dv
\end{gather}

\subsection{Macroscopic constitutive tangents}

The updates of the macrolevel stress and electric displacement are performed by the following incremental relations
\begin{gather}
\Delta \boldsymbol{\sigma}=^{4}\boldsymbol{C}_{\sigma \varepsilon }:\Delta \boldsymbol{\varepsilon}+^{5}\boldsymbol{C}_{\sigma g  } \vdots \Delta \boldsymbol{g}+^{3}\boldsymbol{C}_{\sigma E  }\cdot \Delta (-\boldsymbol{E})\\
\Delta \boldsymbol{\tau }=^{5}\boldsymbol{C}_{\tau \varepsilon }:\Delta \boldsymbol{\varepsilon}+ ^{6}\boldsymbol{C}_{\tau g }\vdots \Delta \boldsymbol{g}+ ^{4}\boldsymbol{C}_{\tau E  }\cdot \Delta (-\boldsymbol{E})\\
\Delta \boldsymbol{D }=^{3}\boldsymbol{C}_{D \varepsilon }:\Delta \boldsymbol{\varepsilon}+ ^{4}\boldsymbol{C}_{D g  } \vdots \Delta \boldsymbol{g}+^{2}\boldsymbol{C}_{D E  }\cdot \Delta (-\boldsymbol{E})
\end{gather}

\begin{figure}[htbp] 
  \centering
  \includegraphics[scale=0.6]{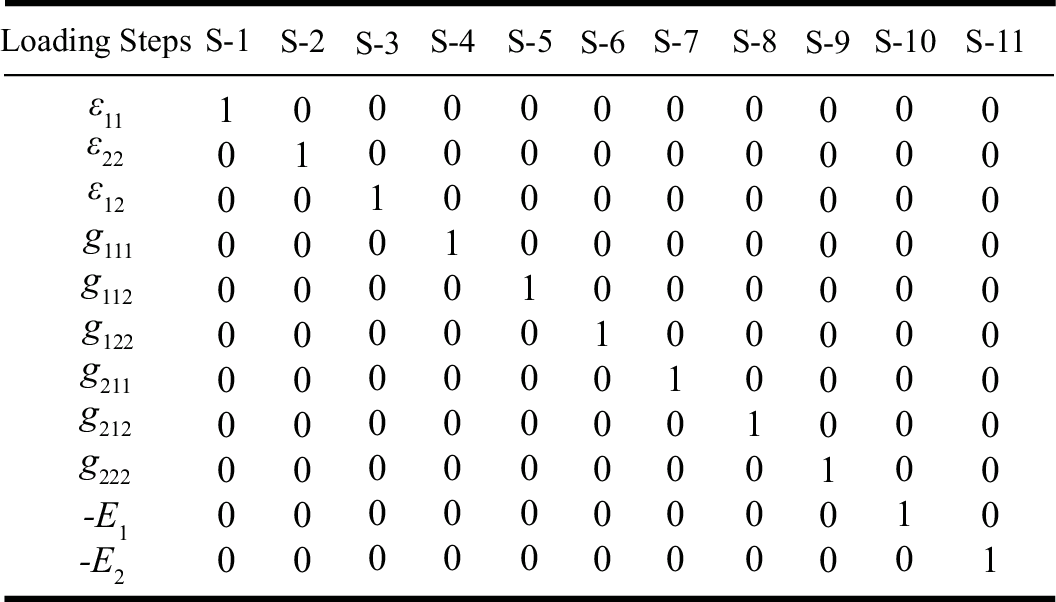}
  \caption{Loading steps for macroscopic constitutive tangents calculation}
  \label{fig.1}
\end{figure}

The tangent stiffness matrices can be calculated with the perturbation analysis by solving linear equations for RVE, associated with incremental changes of macroscopic strains, strain gradients, and electric fields \cite{kaczmarczyk2008scale}. The loading steps for macroscopic constitutive tangents are shown in Fig. \ref{fig.1}. The so-called two-index notation is introduced to represent the homogenized constitutive tangents. 

\section{Numerical implementation}

\subsection{$C^1$ finite elements}

\begin{figure}[htb]
  \centering
  \begin{minipage}[t]{0.48\textwidth}
  \centering
  \includegraphics[scale=0.3]{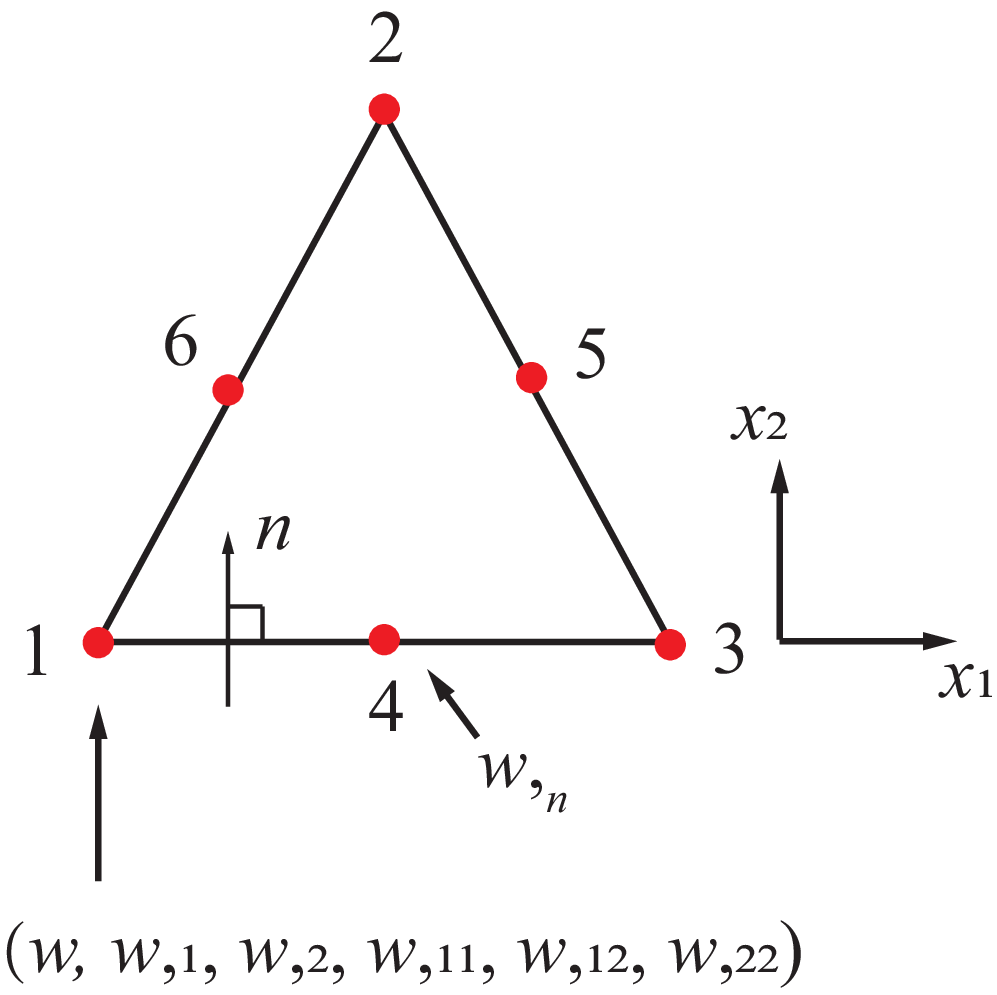}
  \caption{Argyris  triangle  }
  \label{Argyris}
  \end{minipage}
  \begin{minipage}[t]{0.48\textwidth}
  \centering
  \includegraphics[scale=0.30]{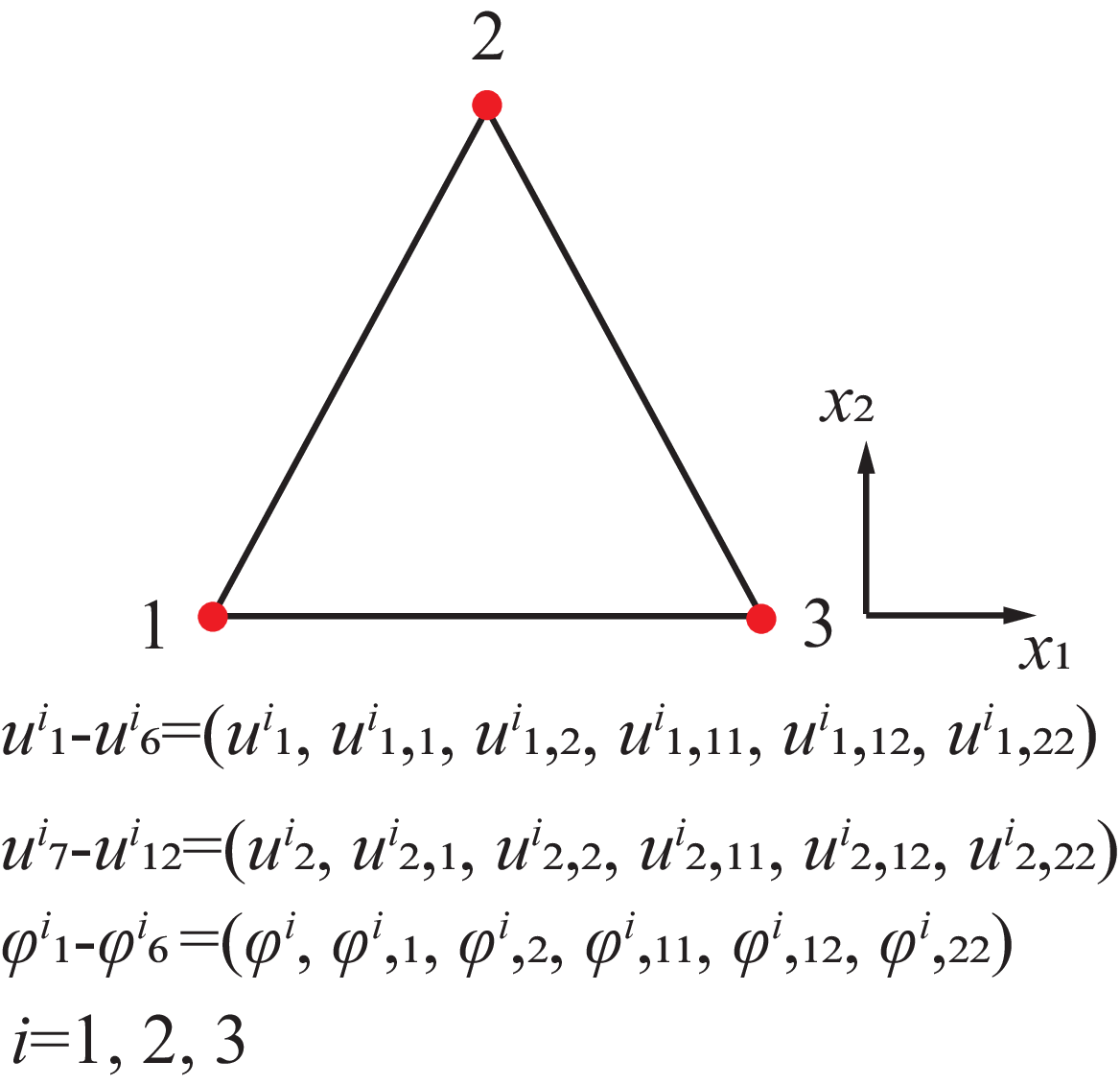}
  \caption{Bell triangle  }
  \label{c1 element}
  \end{minipage}
\end{figure}

 $C^1$ finite elements, that were initially developed for modeling plate structures, are well suited to second-order homogenization of flexoelectricity. It is easier to impose the periodic boundary conditions due to the presence of displacement derivatives in the nodal degrees of freedom. A review and catalog of the $C^1$ finite elements can be found in \cite{hrabok1984review}. The Argyris triangle \cite{argyris1968tuba}  is a well-known $h$-version $C^1$ conforming element. The degrees of freedom include all second-order derivatives on vertices and normal derivatives at mid-side nodes (Fig. \ref{Argyris}), with complete fifth-order polynomial interpolation functions. However, the normal derivatives at midsides are undesirable and are eliminated by interpolating the normal derivative at the mid-side nodes from the derivatives at the corner nodes, resulting in the $C^1$ Bell triangle element \cite{bell1969refined}. The displacement field still varies as a complete quintic inside the element, and the normal derivative along the edges of the element is constrained to be cubic \cite{zervos2001finite}. The governing equation of flexoelectricity is a fourth-order partial differential equation, so the Bell triangle element satisfies the requirements of completeness and compatibility. The shape functions for the Bell triangle were given by Dasgupta and Sengupta \cite{dasgupta1990higher}.  The utility of the Bell triangle was proven in a series of studies on gradient elasticity problems \cite{zervos2001finite, zervos2009two, manzari2013c1}. 

The Bell triangle is applied here as shown in Fig. \ref{c1 element}. Each node has eighteen degrees of freedom, which are displacement and electric potential and their first and second-order derivatives.

In the absence of $E_{i,j}$, $\bar{r}_{i}$ and $\bar{v}_{i}$, the weak form of formulation for flexoelectricity is written as
\begin{gather}
  \int _{v} (\boldsymbol{\sigma} :\delta \boldsymbol{\varepsilon}+\boldsymbol{\tau} \vdots \delta \boldsymbol{g} -\boldsymbol{D} \cdot \delta \boldsymbol{E}) dv=
  \int_{v}(\boldsymbol{f}^{b}\cdot \delta \boldsymbol{u})dv+
  \int_{\varGamma_{t}}(\bar{\boldsymbol{t}}\cdot \delta \boldsymbol{u})d\varGamma-\int_{v} {\rho ^f}\delta \phi dv-\int_{\varGamma_{D}} \omega \delta \phi d\varGamma
\end{gather}  

The displacement $\boldsymbol{u}$ and electric potential $\phi$ fields, as well as their derivatives, can be approximated according to
\begin{gather}
  \boldsymbol{u}= \boldsymbol{N}_{u}\boldsymbol{v}^e ,\;
  \phi= \boldsymbol{N}_{\phi}\boldsymbol{\varphi}^e\\
  \boldsymbol{\varepsilon}=\boldsymbol{B}_{u}\boldsymbol{v}^e,\;
 \boldsymbol{g } =\boldsymbol{H}_{u}\boldsymbol{v}^e   ,\;
    -\boldsymbol{E}=\boldsymbol{B}_{\phi}\boldsymbol{\varphi}^e  
\end{gather}
where $\boldsymbol{\varepsilon}= (\varepsilon_{11},\varepsilon_{22}, \varepsilon_{12})^T$, $\boldsymbol{g}= (g_{111} , g_{112} , g_{122}, g_{211} , g_{212} , g_{222})^T$, $-\boldsymbol{E}= (-E_{1} , -E_{2})^T$. $\boldsymbol{N}_{u}$, $\boldsymbol{N}_{\phi}$ are the shape function matrix, and $\boldsymbol{v}^e$ and $\boldsymbol{\varphi}^e$ are the nodal degrees of freedom. $\boldsymbol{B}_{u}$, $\boldsymbol{B}_{\phi}$, and $\boldsymbol{H}_{u}$ are matrices containing the gradient and Hessian of the corresponding shape functions.

Substituting Eqs. (45)–(46) into Eq. (44) yields
\begin{gather}
  \begin{bmatrix}
    \boldsymbol{K}_{uu}& \boldsymbol{K}_{u\phi}   \\ 
    \boldsymbol{K}_{\phi u}& \boldsymbol{K}_{\phi \phi}  \\
    \end{bmatrix}
  \begin{bmatrix}
     \boldsymbol{v}   \\ 
     \boldsymbol{\varphi} \\
    \end{bmatrix}=
  \begin{bmatrix}
    \boldsymbol{F}_v   \\ 
    \boldsymbol{F}_\varphi \\
    \end{bmatrix}
  \end{gather}
where
\begin{gather}
  \boldsymbol{K}_{uu}=\sum_{e}\int_{v_e}(\boldsymbol{B}^T_{u}\boldsymbol{C}\boldsymbol{B}_{u}+l^{2} \boldsymbol{H}^T_{u}\boldsymbol{Q}\boldsymbol{H}_{u})dv   \\
  \boldsymbol{K}_{u\phi}=\sum_{e}\int_{v_e}(\boldsymbol{B}^T_{u}\boldsymbol{e}\boldsymbol{B}_{\phi}+\boldsymbol{H}^T_{u}\boldsymbol{f}^T\boldsymbol{B}_{\phi })dv\\
  \boldsymbol{K}_{\phi u}=\sum_{e}\int_{v_e}(\boldsymbol{B}^T_{\phi}\boldsymbol{e}^T\boldsymbol{B}_{u}+\boldsymbol{B}^T_{\phi}\boldsymbol{f}\boldsymbol{H}_{u})dv\\
  \boldsymbol{K}_{\phi \phi}=\sum_{e}\int_{v_e}(\boldsymbol{B}^T_{\phi }(-\boldsymbol{\kappa })\boldsymbol{B}_{\phi})dv\\
  \boldsymbol{F}_{u}=\sum_{e}\int_{v_e}(\boldsymbol{N}^{T}_{u} \boldsymbol{f}^{b})dv+\sum_{e}\int_{\varGamma_t}(\boldsymbol{N}^{T}_{u} \bar{\boldsymbol{t}})d\varGamma \\
  \boldsymbol{F}_{\phi }=-\sum_{e}\int_{v_e}(\boldsymbol{N}^{T}_\phi \rho ^f)dv-\sum_{e}\int_{\varGamma_D}(\boldsymbol{N}^{T}_\phi \omega )d\varGamma
\end{gather}
Moreover, $\boldsymbol{C}$, $\boldsymbol{Q}$, $\boldsymbol{\kappa}$, $\boldsymbol{e}$ and $\boldsymbol{f}$ can be written in matrix format as
\begin{gather}
  \boldsymbol{C}=\begin{bmatrix}
    \lambda +2G & \lambda & 0\\ 
    \lambda &  \lambda +2G & 0 \\ 
     0 & 0  & 4G
    \end{bmatrix}, \;
  \boldsymbol{\kappa} =\begin{bmatrix}
      \kappa_{11}& 0\\ 
     0 & \kappa_{33}
     \end{bmatrix} \\
  \boldsymbol{e}^T =\begin{bmatrix}
      0 & 0 & 2e_{15}\\ 
      e_{31} & e_{33} &0 
      \end{bmatrix}, \;
  \boldsymbol{f} =\begin{bmatrix}
        f_{1} +2 f_{2} &0& f_{1} &0& 2f_{2}&0 \\ 
        0&2 f_{2} &0& f_{1} &0&f_{1} +2 f_{2}   
        \end{bmatrix}    \\
   \boldsymbol{Q}=\begin{bmatrix}
    \lambda +2G & 0& 0& 0& \lambda & 0\\ 
    0&\lambda +3G & 0& G& 0& \lambda \\ 
    0&0&G & 0& G& 0 \\ 
    0&G & 0& G& 0 &0 \\ 
    \lambda & 0& G&0&\lambda +3G & 0 \\ 
     0&\lambda &0&0&0& \lambda +2G 
    \end{bmatrix}   
\end{gather}

\subsection{Periodic boundary conditions}

The twelve displacement-related degrees of freedom per node are expressed as $\boldsymbol{u}_g$ and the six electric potential-related degrees of freedom are denoted as $\boldsymbol{\phi}_g$. $\boldsymbol{u}_g$ and $\boldsymbol{\phi}_g$ of nodes on the RVE boundary can be rewritten into the matrix form as 
\begin{gather}
  \boldsymbol{u}_g=\boldsymbol{W} \boldsymbol{\varepsilon} _M+\boldsymbol{S}\boldsymbol{g} _M+ \boldsymbol{r}_{u}\\
  \boldsymbol{\phi}_g=\boldsymbol{L}(-\boldsymbol{E}_M)+\boldsymbol{r}_{\phi}
\end{gather}
where
\begin{gather}
 \boldsymbol{W}=\begin{bmatrix}
    x_1&0&x_2\\1&0&0\\ 0&0&1\\ 0&0&0\\0&0&0\\0&0&0\\
    0&x_2&x_1\\ 0&0&1\\ 0&1&0\\ 0&0&0\\0&0&0\\0&0&0
   \end{bmatrix},    \;
  \boldsymbol{S}=\begin{bmatrix}
    \frac{1}{2}x^2_1&x_1x_2&\frac{1}{2}x^2_2 &0&0&0 \\
     x_1&x_2&0&0&0&0\\ 0&x_1&x_2&0&0&0\\ 1&0&0 &0&0&0\\0&1&0&0&0&0\\0&0&1&0&0&0\\
    0&0&0&\frac{1}{2}x^2_1&x_1x_2&\frac{1}{2}x^2_2\\ 
    0&0&0&x_1&x_2&0\\ 0&0&0&0&x_1&x_2\\ 0&0&0&1&0&0\\0&0&0&0&1&0\\0&0&0&0&0&1
   \end{bmatrix} ,      \;
   \boldsymbol{L}=\begin{bmatrix}
    x_1&x_2\\1&0\\0&1\\0&0\\0&0\\0&0\\
   \end{bmatrix}  
\end{gather}

Here $x_1$ and $x_2$ denote the coordinates of the node. Using the periodicity according to relations (25)-(27) leads to the equations
\begin{gather}
  \boldsymbol{u}_{gR}-\boldsymbol{u}_{gL} =(\boldsymbol{W}_{R}-\boldsymbol{W}_{L})\boldsymbol{\varepsilon} _M+(\boldsymbol{S}_{R}-\boldsymbol{S}_{L})\boldsymbol{g} _M, \;
  \boldsymbol{u}_{gT}-\boldsymbol{u}_{gB} =(\boldsymbol{W}_{T}-\boldsymbol{W}_{B})\boldsymbol{\varepsilon} _M+(\boldsymbol{S}_{T}-\boldsymbol{S}_{B})\boldsymbol{g} _M\\
  \boldsymbol{\phi}_{gR}-\boldsymbol{\phi}_{gL} = (\boldsymbol{L}_{R}-\boldsymbol{L}_{L})(-\boldsymbol{E} _M), \;
  \boldsymbol{\phi}_{gT}-\boldsymbol{\phi}_{gB} = (\boldsymbol{L}_{T}-\boldsymbol{L}_{B})(-\boldsymbol{E} _M)
\end{gather}

The above periodic condition can be called the gradient generalized periodic condition (PBC) \cite{lesivcar2017two}. There is another boundary condition-gradient displacement periodic condition (DBC), where microfluctuations on the RVE boundaries are suppressed, so the displacement and electric potential field of nodes on the RVE boundary are then written as 
\begin{gather}
  \boldsymbol{u}_g=\boldsymbol{W} \boldsymbol{\varepsilon} _M+\boldsymbol{S}\boldsymbol{g} _M\\
  \boldsymbol{\phi}_g=\boldsymbol{L}(-\boldsymbol{E}_M)
\end{gather}

\subsection{The second-order multiscale solution scheme}

\begin{figure}[htb] 
  \centering
  \includegraphics[scale=0.7]{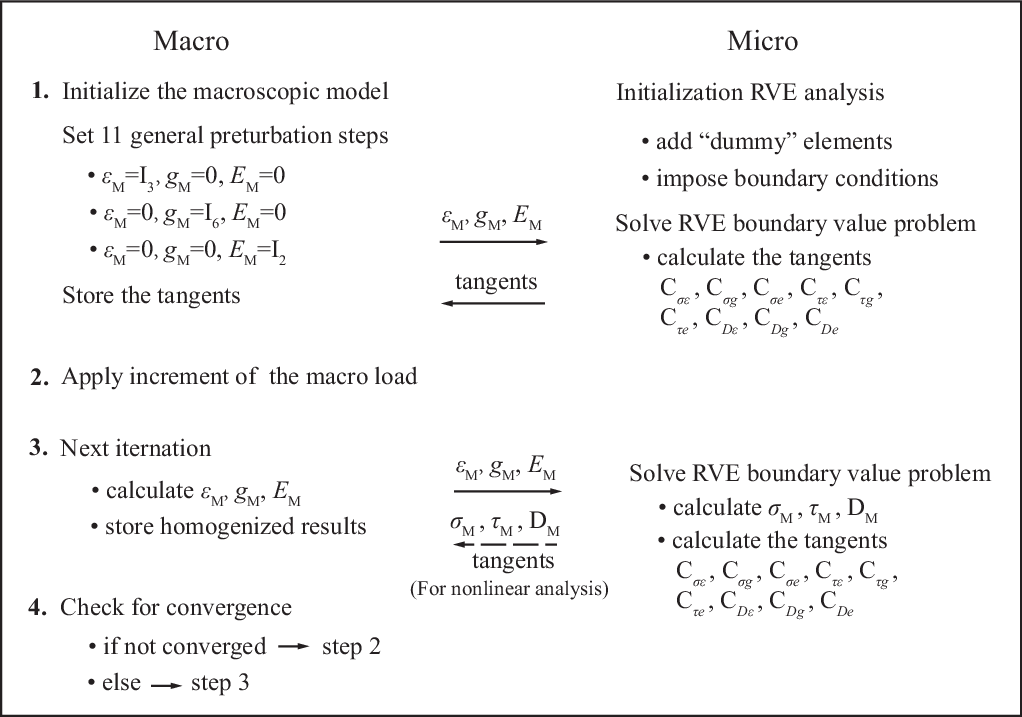}
  \caption{The second-order multiscale solution scheme}
  \label{scheme}
\end{figure}

\begin{figure}[htb] 
  \centering
  \includegraphics[scale=0.7]{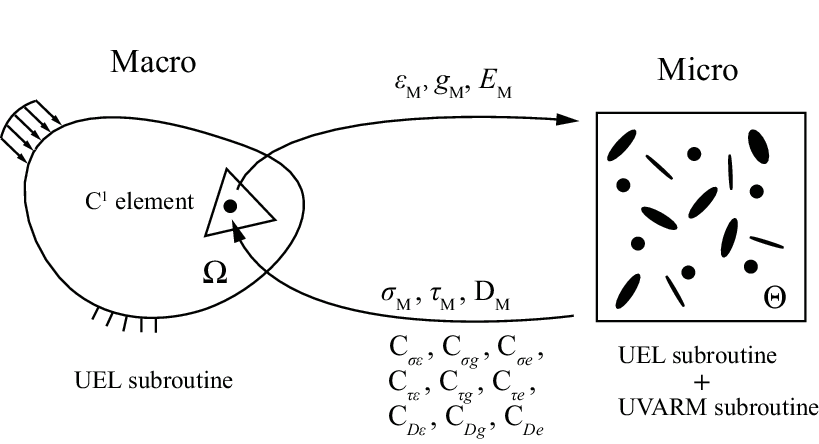}
  \caption{The macro-micro scale transition}
  \label{fig.5}
\end{figure}

The second-order multiscale solution scheme is shown in Fig. \ref{scheme}. The second-order homogenization procedure for flexoelectricity is implemented with ABAQUS. $C^1$ triangular finite element is realized using the user element subroutine UEL, and user subroutine UVARM with ``dummy" elements are used for postprocessing the field results of user-defined elements. The field results are obtained as output in TECPLOT format for visualization. As user subroutine UVARM cannot be used with linear perturbation procedure in ABAQUS, the macroscopic constitutive tangents are calculated with perturbation analysis in the general step. The periodic boundary conditions are imposed with linear constraint equations in ABAQUS. As this study is on flexoelectric structures under small elastic deformation, the macro-micro scale transition (Fig. \ref{fig.5}) is carried out only once. 

\section{ Numerical examples }

\subsection{Material equivalent  parameter analysis of a square plate with holes}

For the sake of comparison, only the representative material equivalent parameters of the flexoelectric microstructure are analyzed, and the equivalent parameters in the following examples are all expressed by the two-index notation. Taking the square plate with holes as an example, the effects of intrinsic length, microscale model size, porosity, and void distribution on the equivalent parameters are investigated. Furthermore, two different boundary conditions (PBC and DBC) are also compared. The material properties of microscale structure are listed in Tab. \ref{Tab.1} based on \cite{deng2017mixed, nanthakumar2017topology}. 
\begin{table}[htbp]
  \centering
  \caption{Material properties ($[\lambda, G]=\rm{GPa}$, $[e]=\rm{C/m^2}$, $[\kappa]=\rm{nC/Vm}$, $[f]=\rm{\mu C/m}$)}
  \begin{tabular}{ccccccccccc}  
    \toprule
    $\lambda$ &  $G$& $e_{31}$& $e_{33}$ & $e_{15}$& $k_{11}$ & $k_{33}$ & $f_{1}$ & $f_{2}$\\
    \midrule
     $179$ & $54$ & $-2.7$ & $3.65$& $21.3$ & $12.5$ & $14.4$ & $1$& $1$\\
    \bottomrule   
  \end{tabular} 
  \label{Tab.1}
\end{table}

\subsubsection{The effects of intrinsic length and microscale model size}

The expression of intrinsic length can be given by the self-consistent estimate \cite{aifantis2000gradient} and its calculation is beyond the scope of this paper. The intrinsic length in flexoelectricity analysis is taken as the model size multiplied by a factor smaller than one \cite{ deng2017mixed, mao2016mixed}, but the value of the factor is not definite.

\begin{figure}[htbp] 
  \centering
  \includegraphics[scale=0.3]{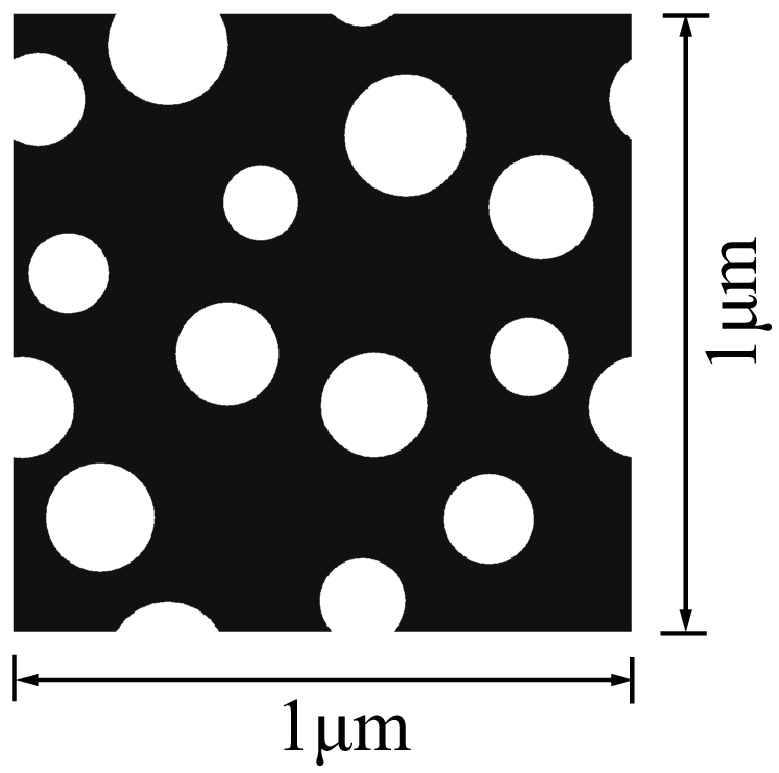}
  \caption{The mircroscale model with holes}
  \label{fig.6}
\end{figure}

\begin{figure}[htb]
  \centering
  \begin{minipage}[t]{0.48\textwidth}
  \centering
  \includegraphics[scale=0.30]{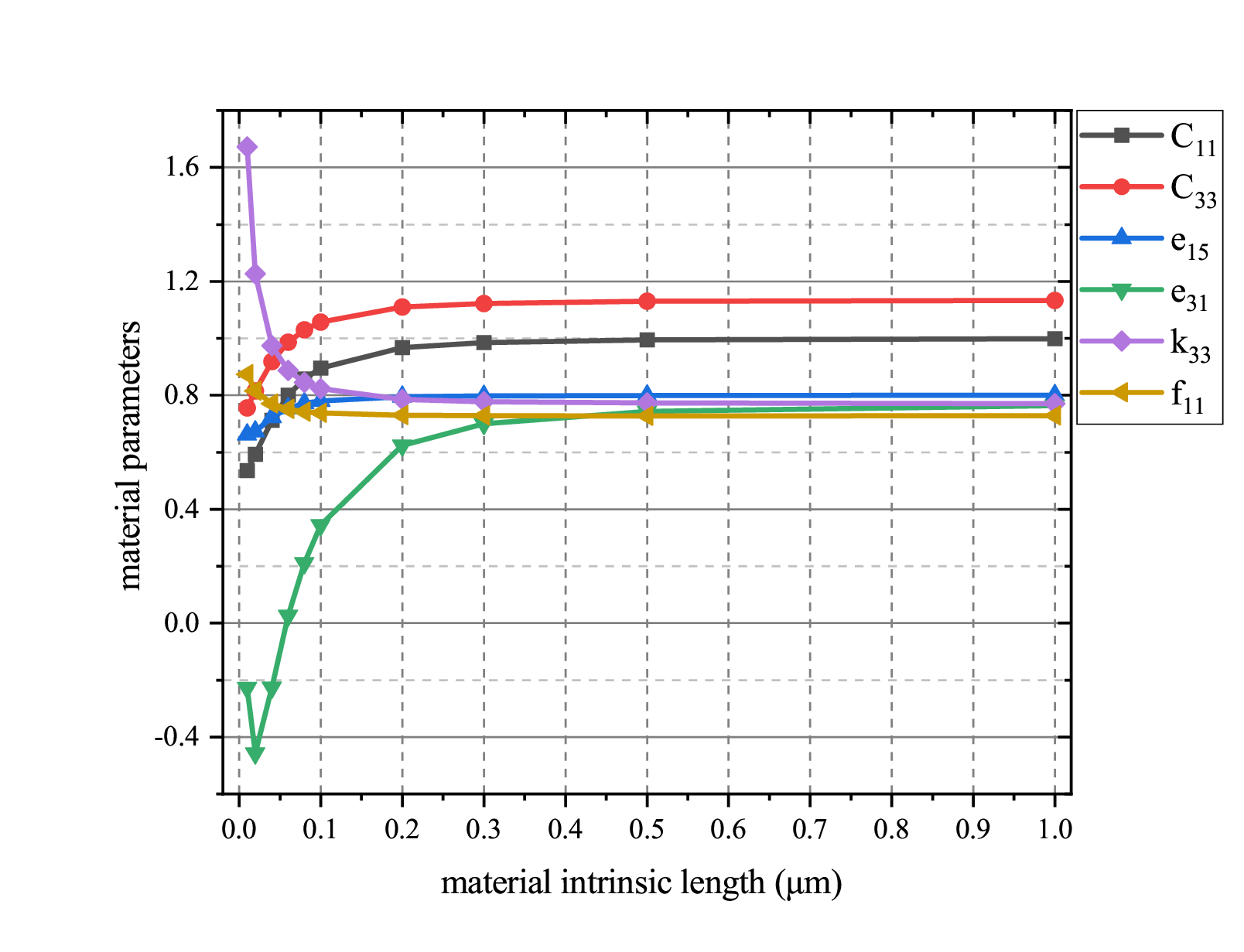}
  \caption{Material parameters of different intrinsic length l }
  \label{fig.7}
  \end{minipage}
  \begin{minipage}[t]{0.48\textwidth}
  \centering
  \includegraphics[scale=0.30]{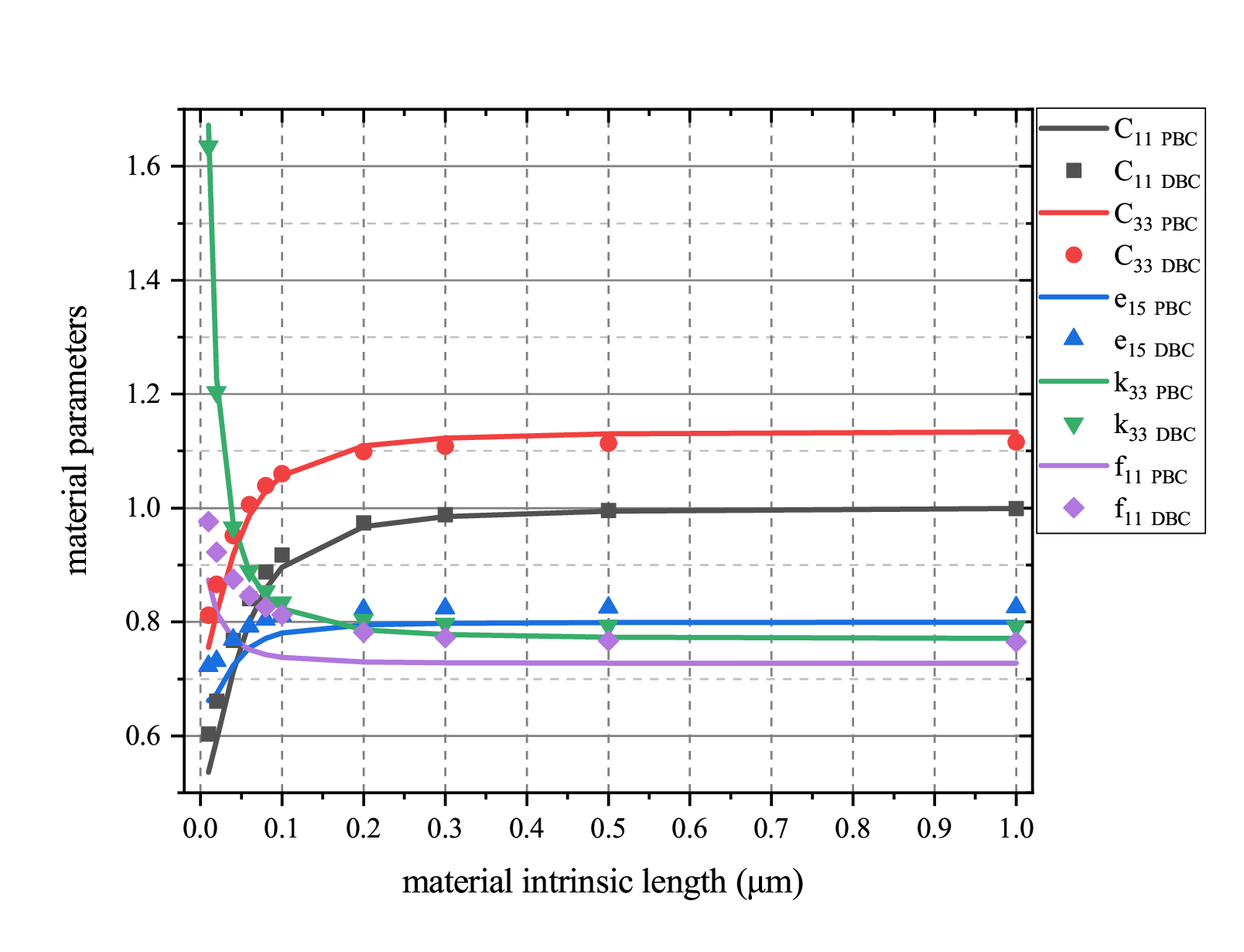}
  \caption{Material parameters of different intrinsic length l with different boundary conditions}
  \label{fig.8}
  \end{minipage}
\end{figure}

The microscale model with a side length of 1 $\mu \rm{m}$ is shown in Fig. \ref{fig.6}. The effects of intrinsic length on the equivalent parameters of microstructure are shown in Fig. \ref{fig.7}. The parameters are normalized by dividing the corresponding material property in Tab. \ref{Tab.1} to aid comparison. One can observe that the dielectric and flexoelectric coefficients decrease and the elastic and piezoelectric coefficients increase with the increase of intrinsic length, while all the coefficients remain unchanged when the intrinsic length is over 0.5 $\mu \rm{m}$. The effects of boundary conditions, i.e. PBC and DBC, on the results are shown in Fig. \ref{fig.8}. The equivalent elastic coefficients remain quite similar irrespective of the considered boundary conditions, while the dielectric, piezoelectric, and flexoelectric coefficients differ noticeably under DBC.


\begin{figure}[htbp]
  \centering
  \begin{minipage}[t]{0.48\textwidth}
  \centering
  \includegraphics[scale=0.3]{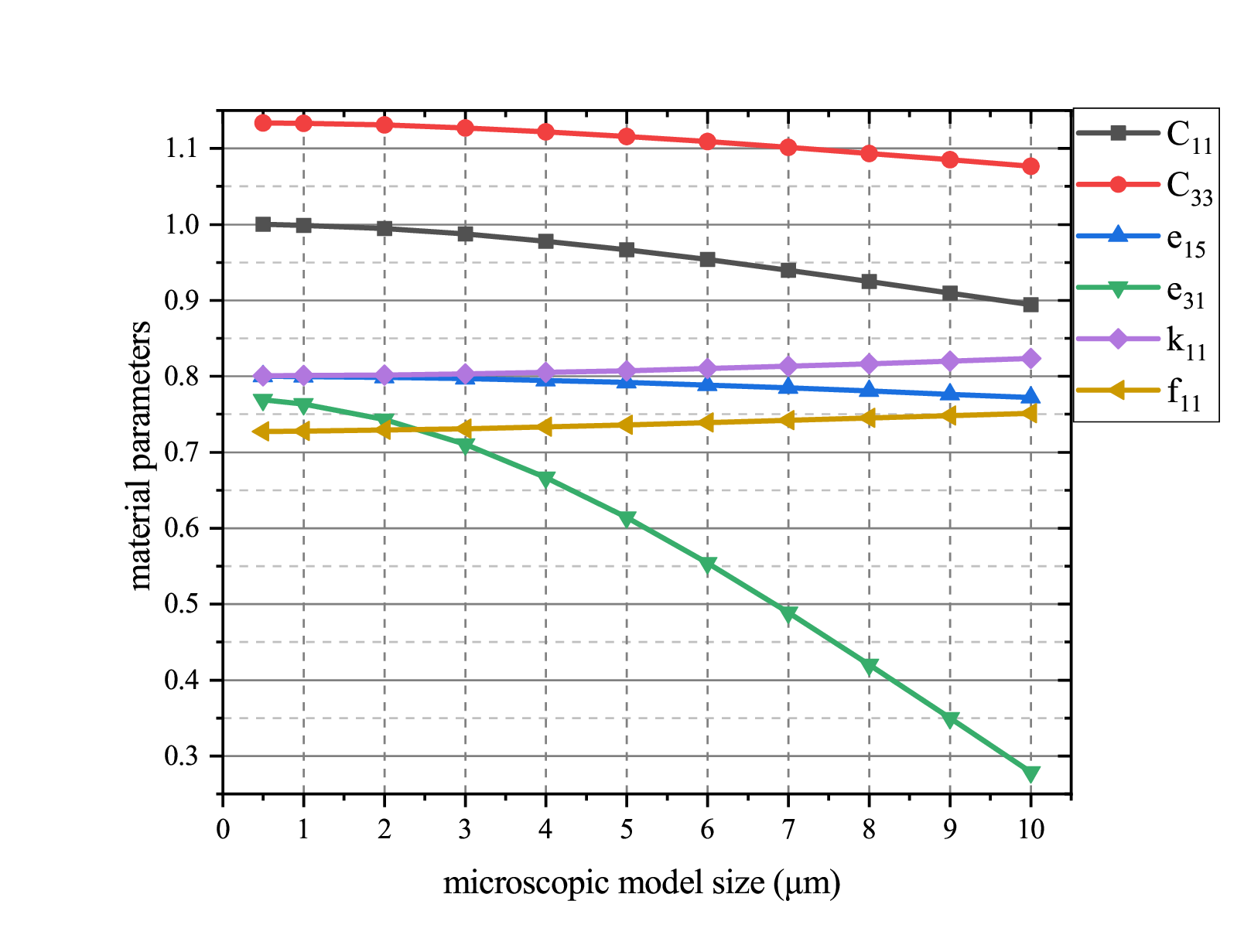}
  \caption{Material parameters of different microscopic model sizes }
  \label{fig.9}
  \end{minipage}
  \begin{minipage}[t]{0.48\textwidth}
  \centering
  \includegraphics[scale=0.3]{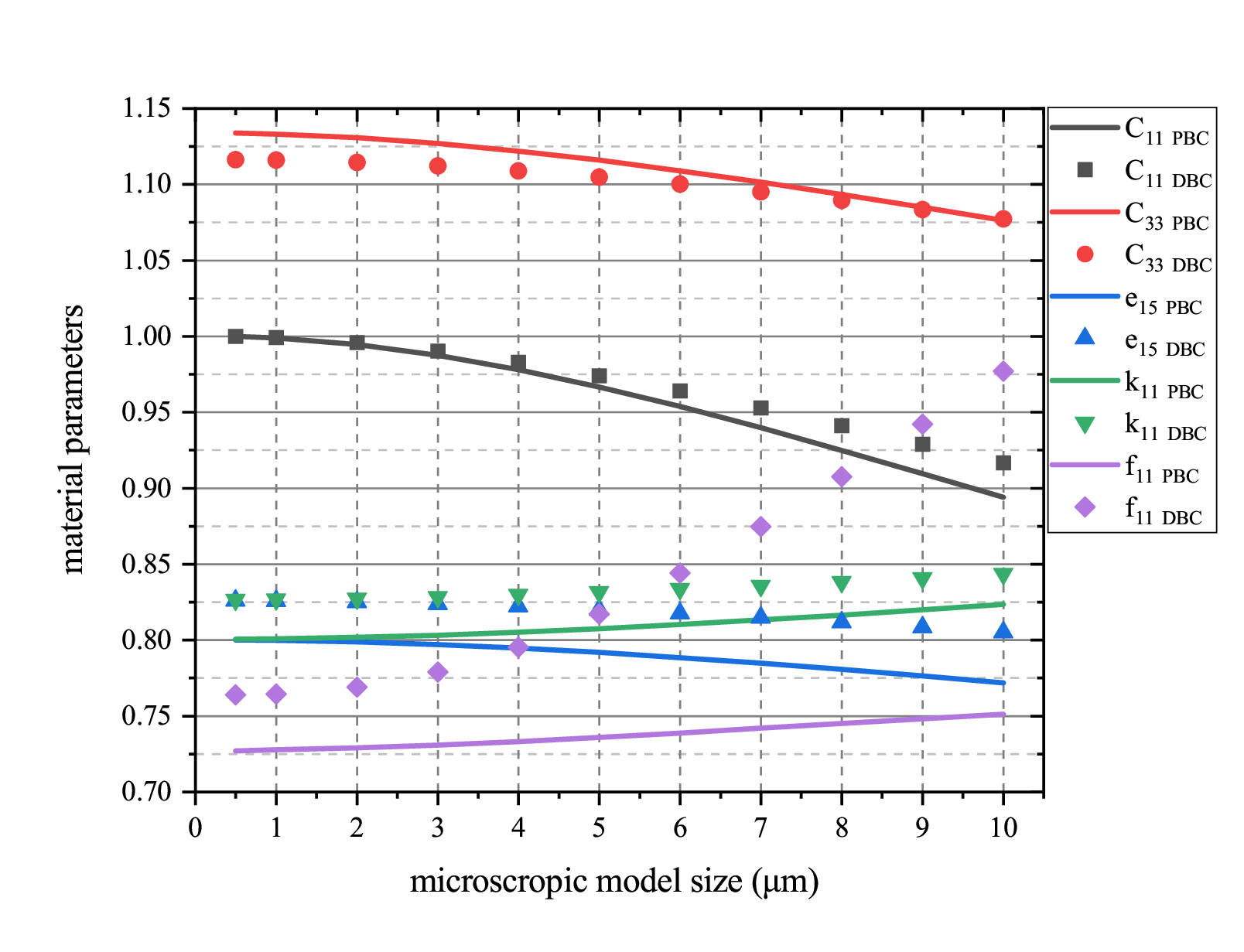}
  \caption{Material parameters of different microscopic model sizes with different boundary conditions}
  \label{fig.10}
  \end{minipage}
\end{figure}

Keeping the intrinsic length constant ($1 \rm{\mu m}$),  the effects of microscale model size are shown in Fig. \ref{fig.9}.  The elastic and piezoelectric coefficients decrease, while the dielectric and flexoelectric coefficients increase with the increase of microscale model size. On the other hand, the dielectric, piezoelectric, and flexoelectric coefficients are also larger under DBC compared to PBC, as in Fig. \ref{fig.10}. 

\subsubsection{The effects of porosity and void distribution}

\begin{figure}[htbp]
  \centering
  \begin{minipage}[t]{0.48\textwidth}
  \centering
  \includegraphics[scale=0.3]{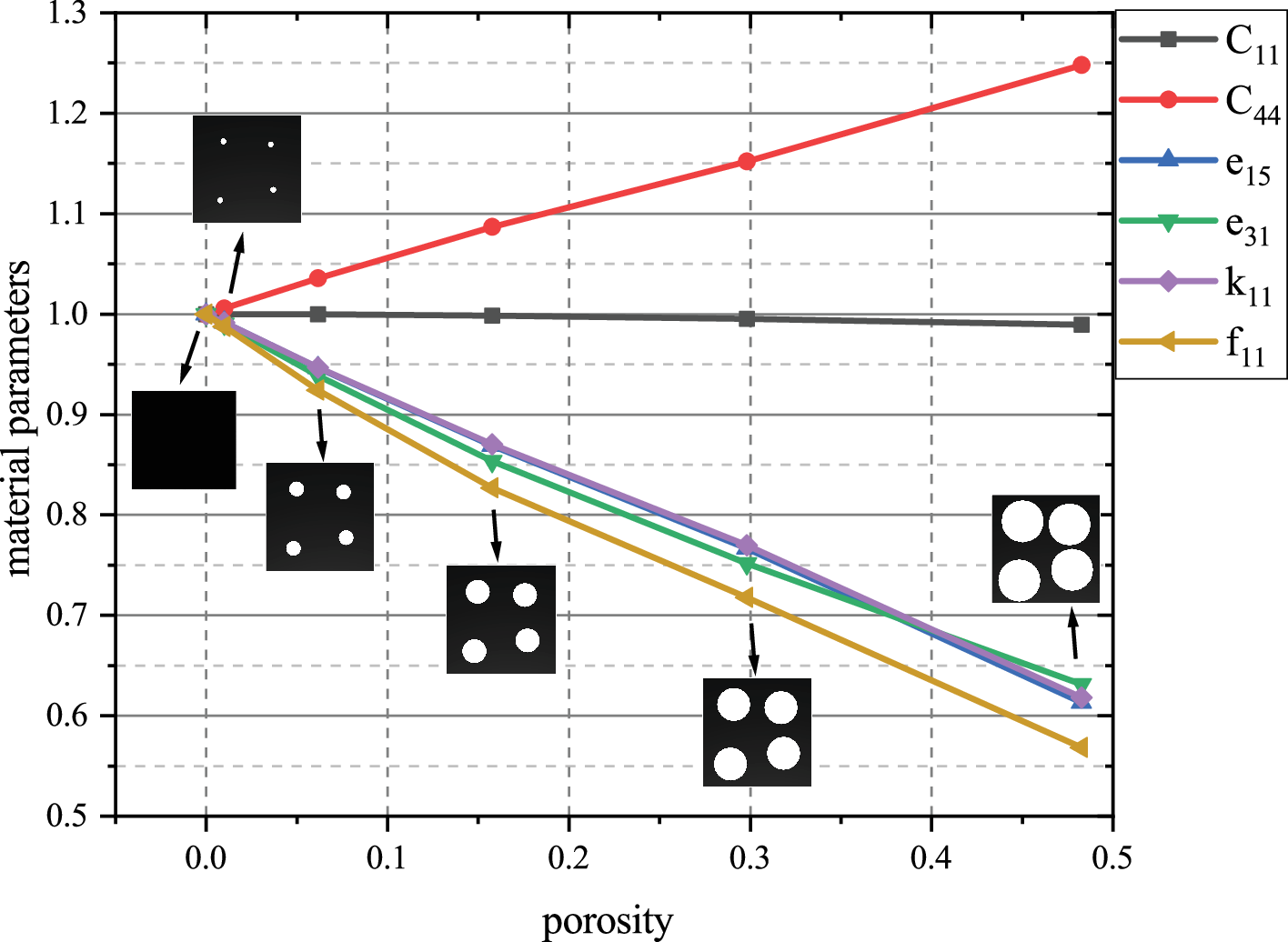}
  \caption{Material parameters of different porosities }
  \label{fig.11}
  \end{minipage}
  \begin{minipage}[t]{0.48\textwidth}
  \centering
  \includegraphics[scale=0.3]{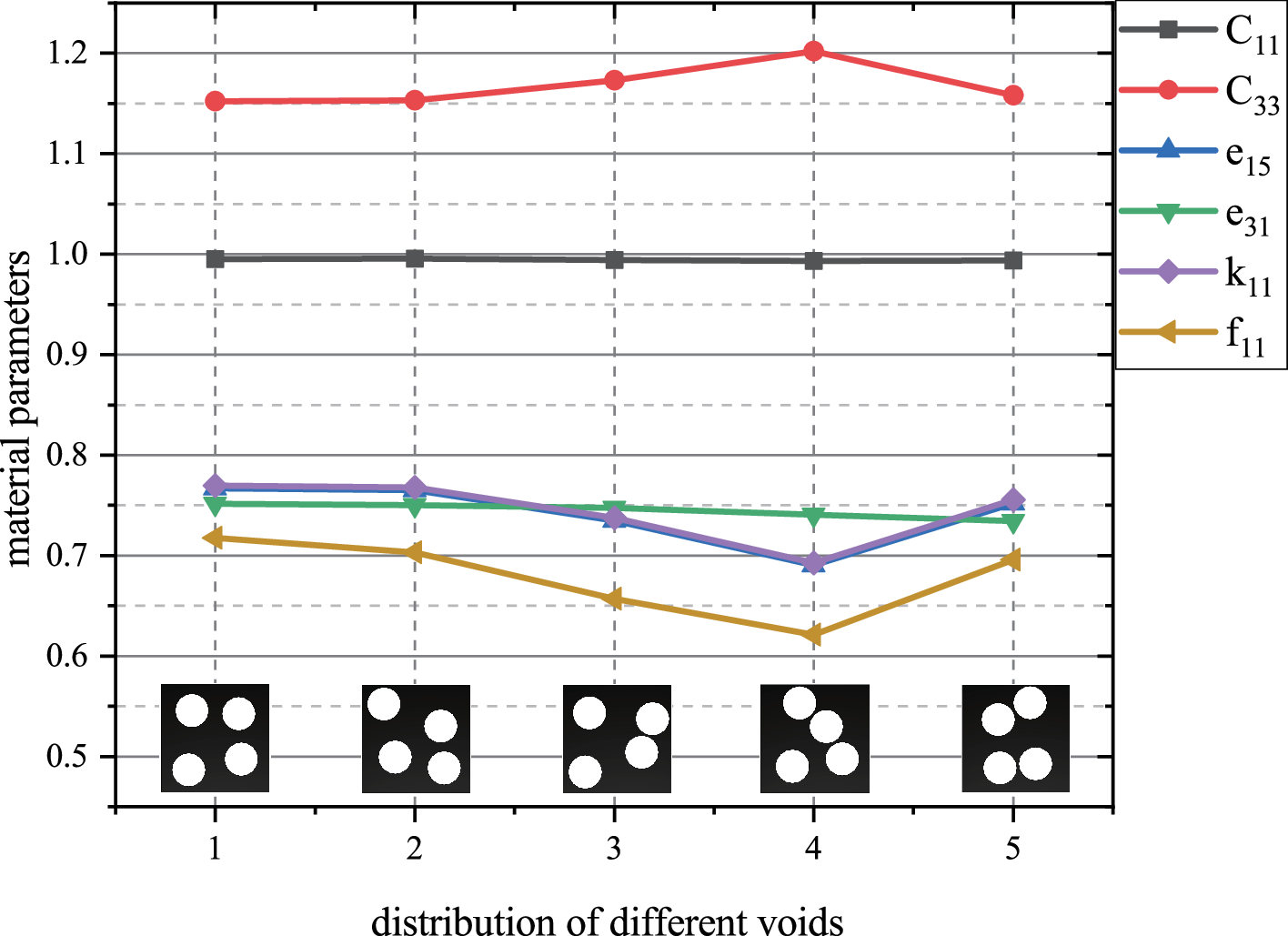}
  \caption{Material parameters of different distribution of holes}
  \label{fig.12}
  \end{minipage}
\end{figure}

The side length of the model and intrinsic length are $1 \ \mu \rm{m}$, and PBC is set as the boundary condition. The variation of the equivalent material parameters is linearly related to the porosity as shown in Fig. \ref{fig.11}. As the porosity increases, the dielectric, piezoelectric and flexoelectric coefficients decrease, and the elastic coefficient $C_{33}$ increases, while the elastic coefficient $C_{11}$ remains almost unchanged. 
It should be noted that when the porosity is zero, the equivalent material parameters are the same as the microscopic material, which is consistent with intuition.
According to Fig. \ref{fig.12}, the equivalent flexoelectric parameter is more obviously impacted by the void distribution than other parameters. Compared with void distribution, the influence of shape and topology of microstructure on the equivalent material parameters will be greater, which needs further study.

\subsection{Creating materials with piezoelectric behaviour without using piezoelectric materials}
Flexoelectricity can be exploited to design piezoelectric effects without using piezoelectric materials by seeking to induce strain gradients \cite{krichen2016flexoelectricity, deng2014flexoelectricity}. One method is to introduce a triangular hole in a sheet, then a net nonzero polarization exists under uniform stretching. It is to be noted that, the average polarization remains zero if the shape of the hole is circular, as shown in Fig. \ref{fig.13}.  In this paper, similar methods are applied in the microscale model to create equivalent piezoelectric coefficients in the absence of piezoelectric materials. 

\subsubsection{Two-scale linear analysis of a uniaxial tensioned square plate}

\begin{figure}[htb]
  \centering
  \begin{minipage}[t]{0.48\textwidth}
  \centering
  \includegraphics[scale=0.6]{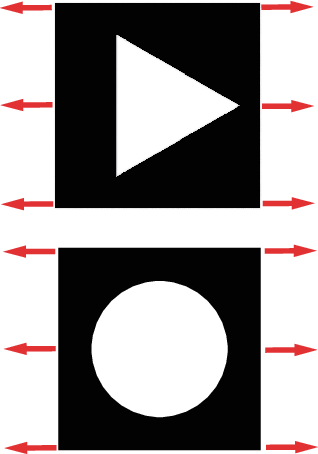}
  \caption{Piezoelectric nanocomposites without using piezoelectric materials \cite{krichen2016flexoelectricity}}
  \label{fig.13}
  \end{minipage}
  \begin{minipage}[t]{0.48\textwidth}
  \centering
  \includegraphics[scale=0.4]{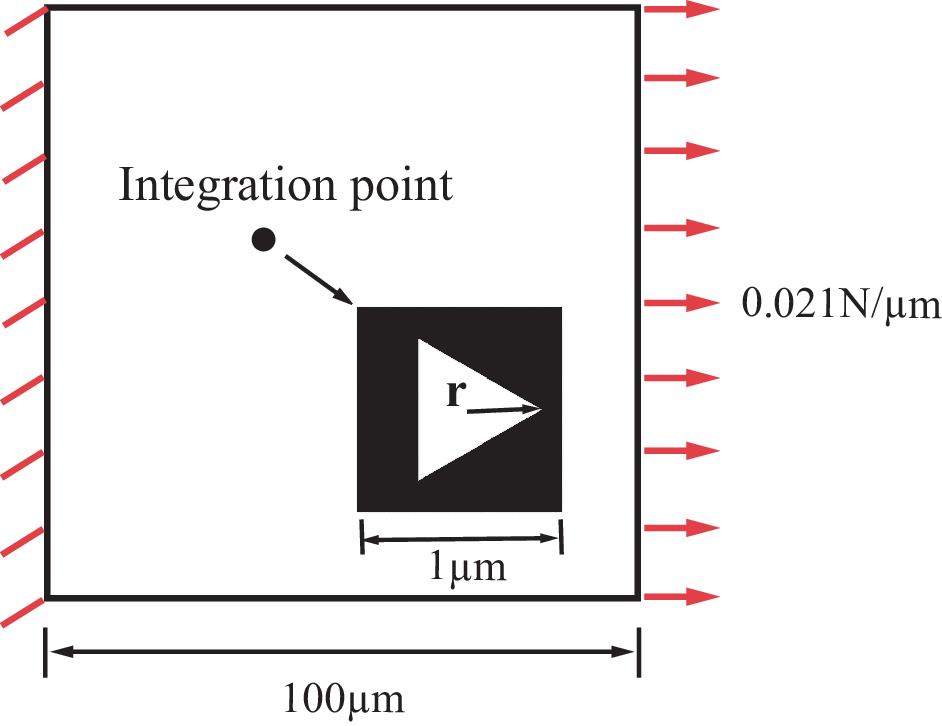}
  \caption{Macroscopic and microscopic models}
  \label{fig.14}
  \end{minipage}
\end{figure}

The macroscale model is a square plate subjected to uniaxial tension, and the microscale model is a square with a triangular hole in the center (Fig. \ref{fig.14}). The material properties are similar to the ones used in the previous examples as listed in Tab. \ref{Tab.1}, except that the values of piezoelectric coefficients are considered to be zero. The intrinsic length is set to be 1 $ \rm{ \mu m}$.

\begin{figure}[htbp]
  \centering
  \begin{minipage}[t]{0.48\textwidth}
  \centering
  \includegraphics[scale=0.3]{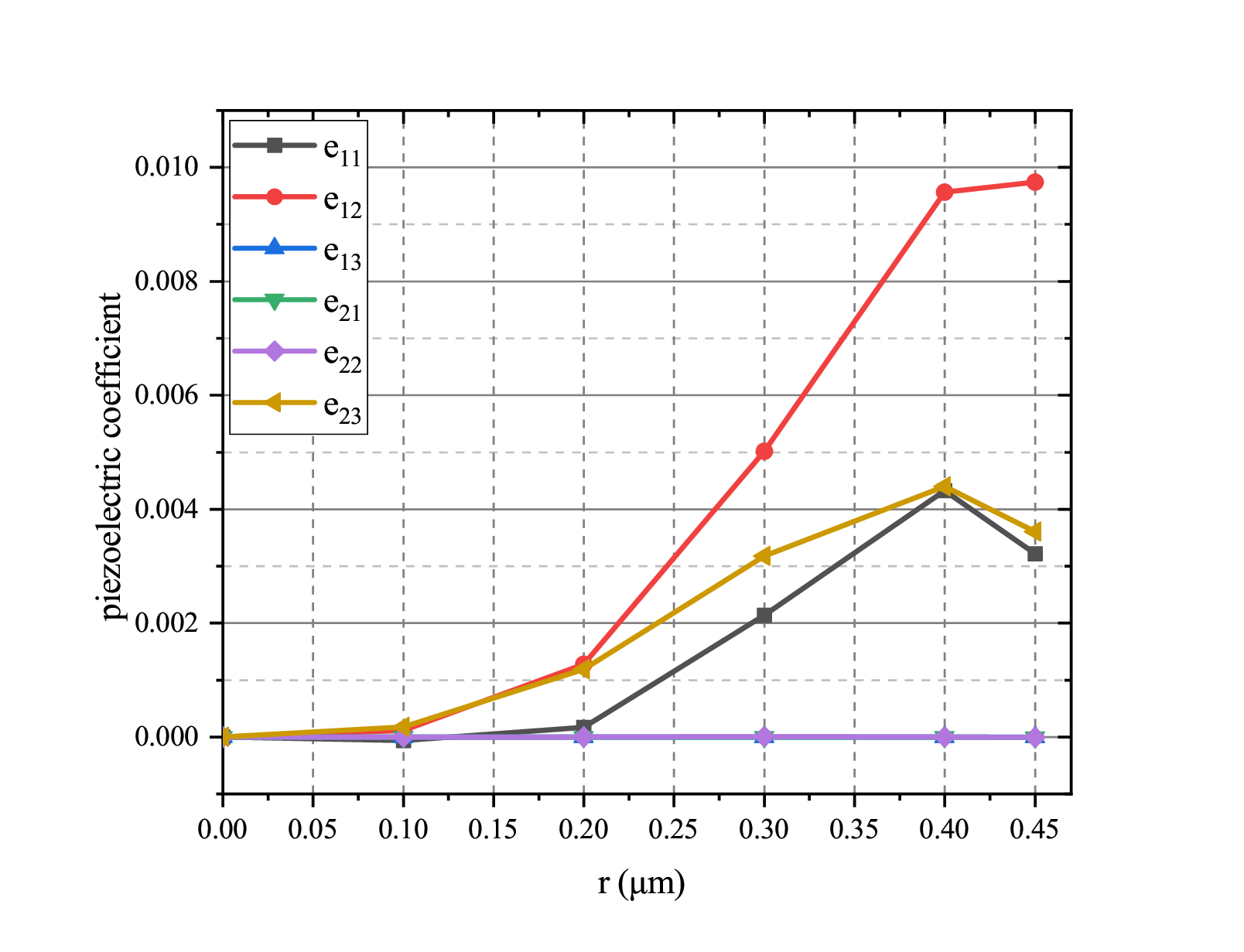}
  \caption{Equivalent piezoelectric coefficients ($\rm{C/m^2}$) of different hole sizes}
  \label{fig.15}
  \end{minipage}
  \begin{minipage}[t]{0.48\textwidth}
  \centering
  \includegraphics[scale=0.4]{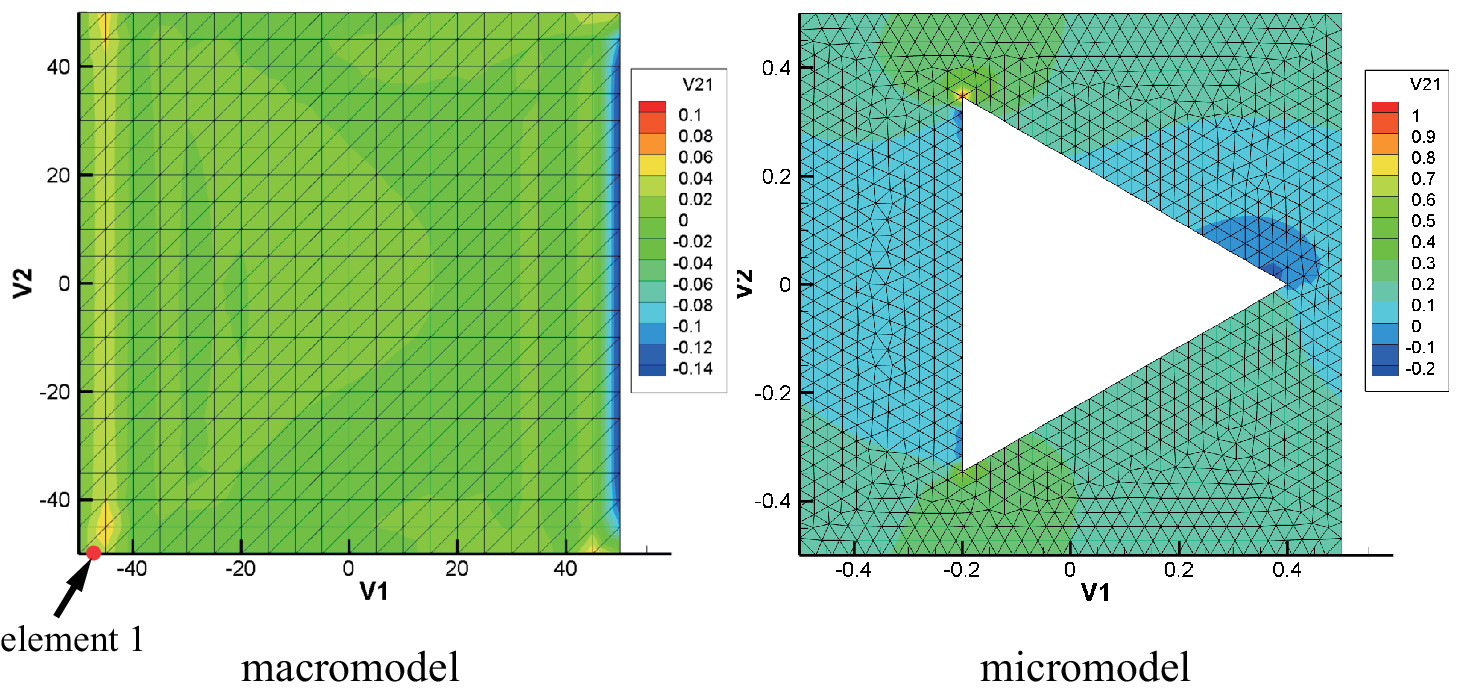}
  \caption{$D_1$ of the macroscale and microscale model ($\rm{C/m^2}$)}
  \label{Macromicro}
  \end{minipage}
\end{figure}

Fig. \ref{fig.15} gives piezoelectric coefficients with different hole sizes measured by the radius $r$ of the inscribed circle of a triangle. The piezoelectric coefficients $e_{11}$, $e_{12}$, $e_{23}$ generally increase as $r$ increases, while $e_{13}$, $e_{21}$, $e_{22}$ are almost zero. It is noteworthy that the equivalent piezoelectric coefficients are calculated to be zero when the hole is a circle, regardless of the hole size, which is consistent with the discussion in \cite{krichen2016flexoelectricity}.

The electric displacement field $D_1$ of the macro model with $r=0.4  \ \mu  \rm{m}$ in the microstructure is shown in Fig. \ref{Macromicro}. The electric displacement field $D_1$ of the microstructure which corresponds to one integration point of element 1 in the macro model is also given in the figure. The calculation procedure has been described in Fig. \ref{scheme}, and the macro-micro scale transition is carried out only once because of the linear elastic analysis.

\subsubsection{Creating piezoelectric materials with flexoelectric composites}
\begin{figure}[htbp] 
  \centering
  \includegraphics[scale=0.25]{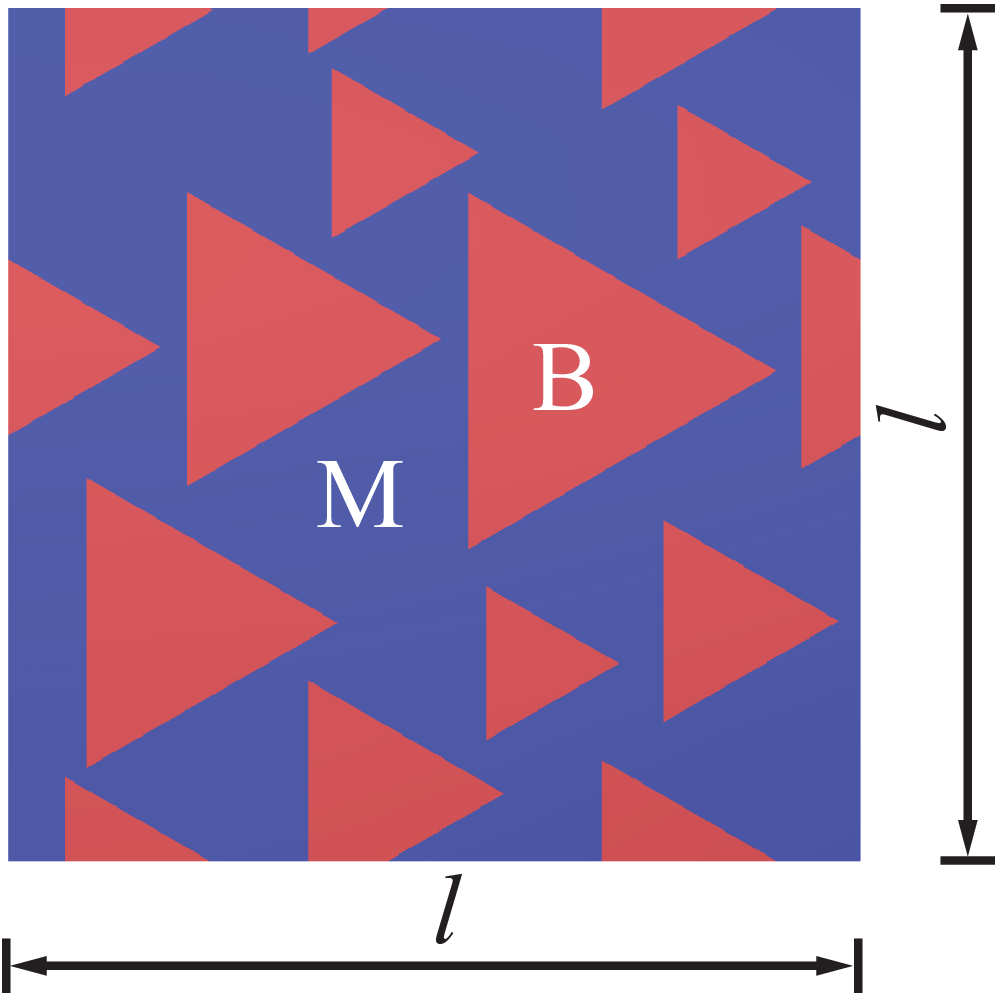}
  \caption{Flexoelectric composites}
  \label{BM}
\end{figure}

\begin{figure}[htb]
  \centering
  \begin{minipage}[t]{0.48\textwidth}
  \centering
  \includegraphics[scale=0.3]{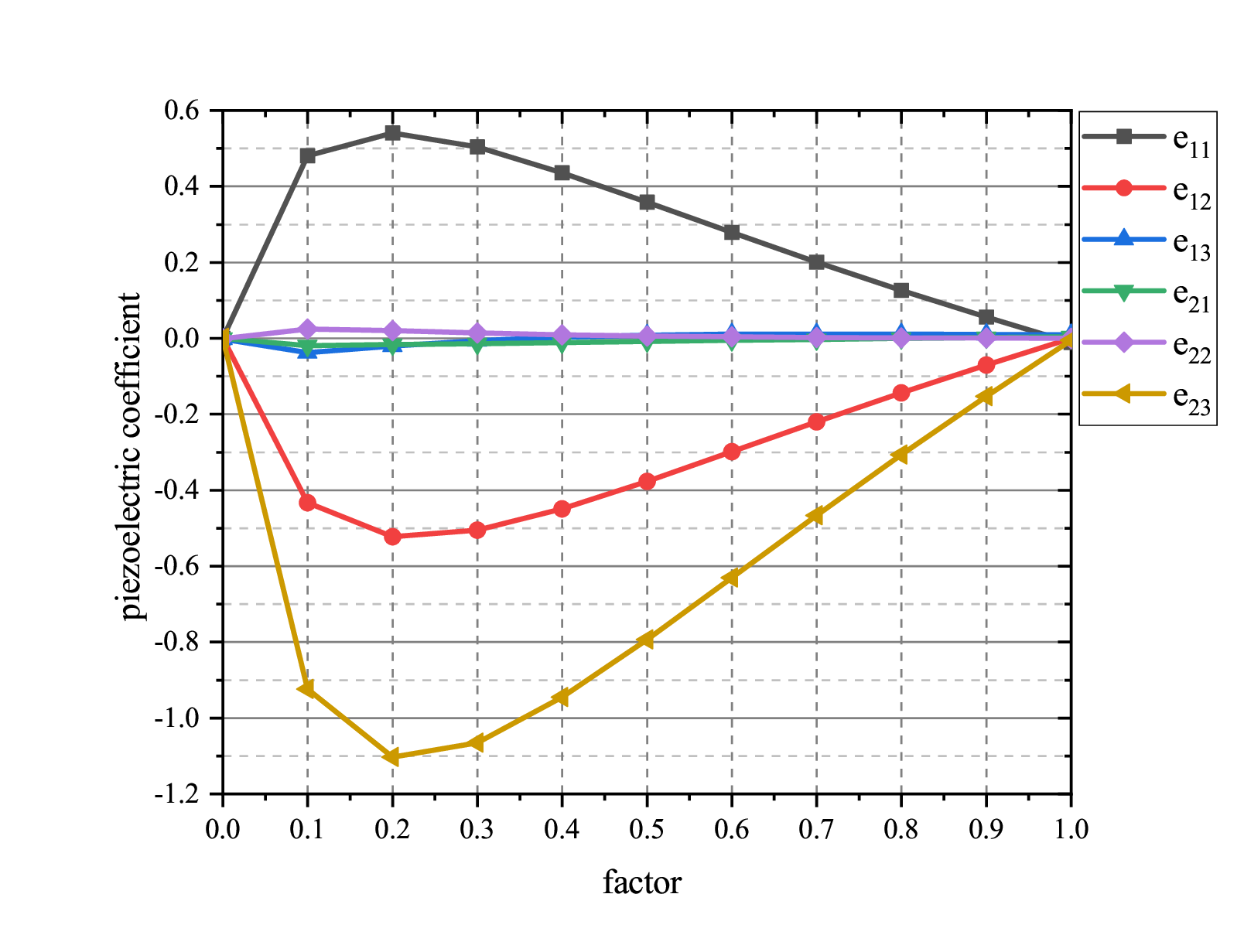}
  \caption{Homogenized piezoelectric coefficients ($\rm{C/m^2}$) of different factors}
  \label{factor}
  \end{minipage}
  \begin{minipage}[t]{0.48\textwidth}
  \centering
  \includegraphics[scale=0.3]{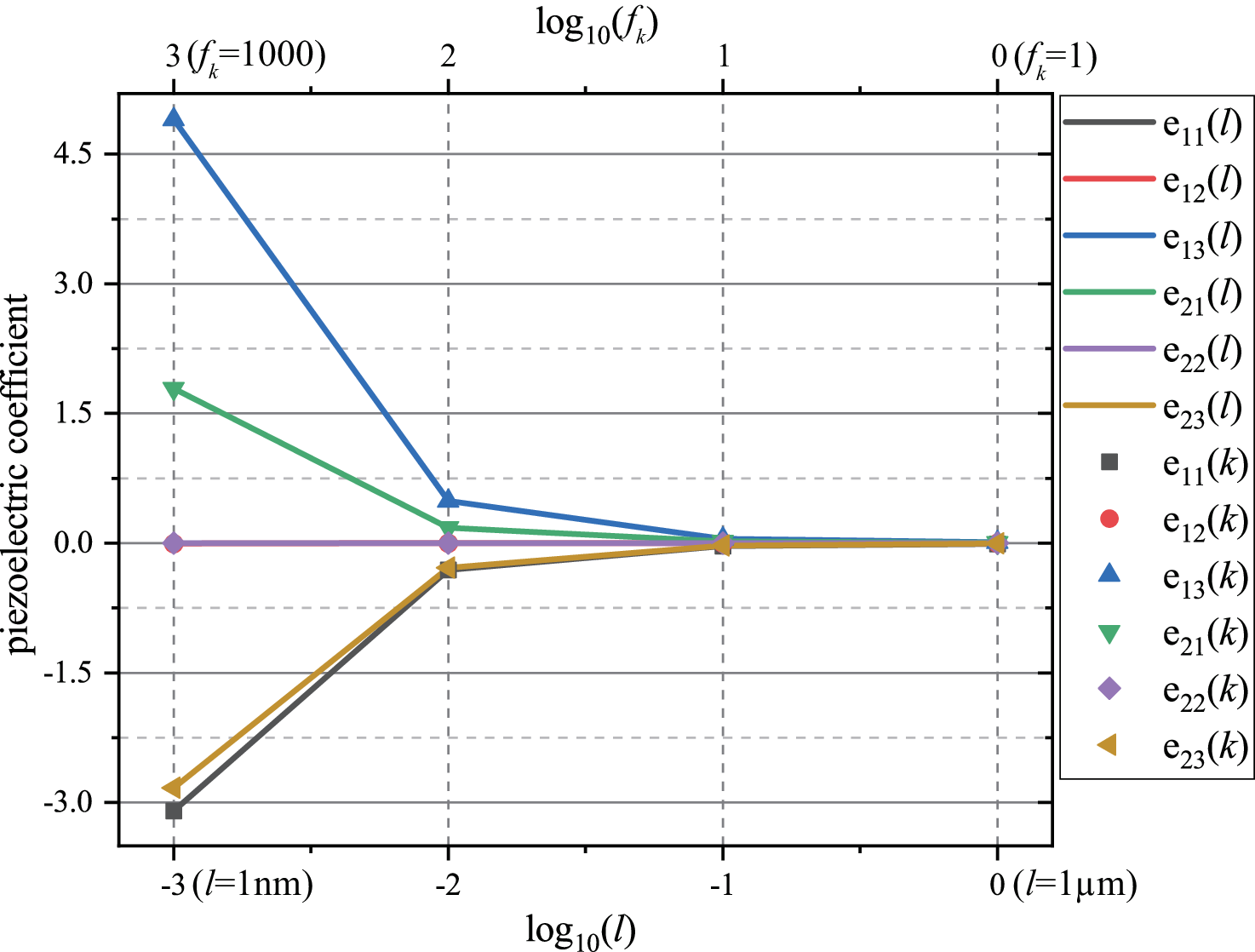}
  \caption{Homogenized piezoelectric coefficients ($\rm{C/m^2}$) of different micromodel sizes $l$ and factors $f_k$}
  \label{size}
  \end{minipage}
\end{figure}

Similarly, we can design piezoelectric materials with flexoelectric composites by distributing flexoelectric inclusions (B) in the dielectric matrix (M), as shown in Fig. {\ref{BM}}. The material properties of the inclusions are listed in Tab. \ref{Tab.1}, the side length is $1 \ \mu \rm{m}$, and the intrinsic length is ignored. The dielectric coefficients of the matrix are the same as the inclusions, but the flexoelectric coefficients are zero. To induce the strain gradients, the elastic tensor of the matrix is adjusted and set to be the inclusions' elastic tensor multiplied by a factor less than 1. As the factor varies between 0 and 1, the homogenized piezoelectric coefficients vary as shown in Fig. \ref{factor}. when the factor is 1 or 0, the homogenized piezoelectric coefficients are nearly zero because the strain gradient is too small (factor=1) or the connection between the inclusions is extremely weak (factor=0). The ideal factor is around 0.2. 

Moreover, the size-dependent behaviour of flexoelectricity can be used to enhance the homogenized piezoelectric coefficients. By reducing the size of the microscale composite structure, from the micron-scale to the nanoscale, the homogenized piezoelectric coefficients increase significantly even though the elastic tensor factor is 1 (Fig. \ref{size}). Interestingly, we can get similar results by increasing the flexoelectric coefficients of the inclusions without shrinking the microstructure size. The same homogenized piezoelectric coefficients can be obtained simply by increasing the flexoelectric coefficients of the inclusions by a factor of $f_k$ (Fig. \ref{size}), instead of scaling the microstructure down when regardless of the intrinsic length, which are consistent with the dimensional analysis. The results in Fig. \ref{size} illustrate the size-dependent behaviour of flexoelectricity in another way. 

\textit{Remark: the self-consistency of the proposed method is reflected in four aspects: 1. The equivalent material parameters are the same as the microscopic material when the porosity is zero (Fig. 11); 2. The equivalent piezoelectric coefficients are calculated to be zero when the hole is a circle; 3. The homogenized piezoelectric coefficients are nearly zero when the elastic tensor factor of the matrix is 1 or 0 (Fig. 18); 4. The same homogenized piezoelectric coefficients can be obtained by increasing the flexoelectric coefficients of the inclusions or shrinking the microstructure size (Fig. 19).}

\section{Conclusion}

The computational second-order homogenization strategy of flexoelectric composites has been presented. The framework was implemented with $C^1$ triangular finite elements based on the secondary development of ABAQUS via user subroutines UEL and UVARM. Compared to the asymptotic expansion method, the homogenized electromechanical properties can be computed using perturbation analysis without complicated mathematical derivations. The influence of intrinsic length, micromodel size, and the two boundary conditions on the effective coefficients of flexoelectric microstructure are discussed. The influence of porosity and void distribution in the microstructure on the equivalent piezoelectric coefficients are also studied. The homogenized piezoelectric coefficients by introducing triangular voids or flexoelectric inclusions in the microstructure are determined. In this work, we restrict the computational second-order homogenization scheme to the linear 2D case, however, considering nonlinearity and extension to 3D problems are straightforward.

\section*{Acknowledgements}
The authors acknowledge the support from ERC Starting Grant (802205).
\section*{References}
\bibliographystyle{elsarticle-num-names}
\bibliography{lite_flexo.bib}

\begin{thebibliography}{45}
\providecommand{\natexlab}[1]{#1}
\providecommand{\url}[1]{\texttt{#1}}
\providecommand{\urlprefix}{URL }
\expandafter\ifx\csname urlstyle\endcsname\relax
  \providecommand{\doi}[1]{doi:\discretionary{}{}{}#1}\else
  \providecommand{\doi}[1]{doi:\discretionary{}{}{}\begingroup
  \urlstyle{rm}\url{#1}\endgroup}\fi
\providecommand{\bibinfo}[2]{#2}

\bibitem[{Majdoub et~al.(2008)Majdoub, Sharma, and Cagin}]{majdoub2008enhanced}
\bibinfo{author}{M.~Majdoub}, \bibinfo{author}{P.~Sharma},
  \bibinfo{author}{T.~Cagin}, \bibinfo{title}{Enhanced size-dependent
  piezoelectricity and elasticity in nanostructures due to the flexoelectric
  effect}, \bibinfo{journal}{Physical Review B}
  \bibinfo{volume}{77}~(\bibinfo{number}{12}) (\bibinfo{year}{2008})
  \bibinfo{pages}{125424}.

\bibitem[{Nguyen et~al.(2013)Nguyen, Mao, Yeh, Purohit, and
  McAlpine}]{nguyen2013nanoscale}
\bibinfo{author}{T.~D. Nguyen}, \bibinfo{author}{S.~Mao},
  \bibinfo{author}{Y.-W. Yeh}, \bibinfo{author}{P.~K. Purohit},
  \bibinfo{author}{M.~C. McAlpine}, \bibinfo{title}{Nanoscale
  flexoelectricity}, \bibinfo{journal}{Advanced Materials}
  \bibinfo{volume}{25}~(\bibinfo{number}{7}) (\bibinfo{year}{2013})
  \bibinfo{pages}{946--974}.

\bibitem[{Zubko et~al.(2013)Zubko, Catalan, and
  Tagantsev}]{zubko2013flexoelectric}
\bibinfo{author}{P.~Zubko}, \bibinfo{author}{G.~Catalan},
  \bibinfo{author}{A.~K. Tagantsev}, \bibinfo{title}{Flexoelectric effect in
  solids}, \bibinfo{journal}{Annual Review of Materials Research}
  \bibinfo{volume}{43}.

\bibitem[{Yudin and Tagantsev(2013)}]{yudin2013fundamentals}
\bibinfo{author}{P.~Yudin}, \bibinfo{author}{A.~Tagantsev},
  \bibinfo{title}{Fundamentals of flexoelectricity in solids},
  \bibinfo{journal}{Nanotechnology} \bibinfo{volume}{24}~(\bibinfo{number}{43})
  (\bibinfo{year}{2013}) \bibinfo{pages}{432001}.

\bibitem[{Wang et~al.(2019)Wang, Gu, Zhang, and
  Chen}]{wang2019flexoelectricity}
\bibinfo{author}{B.~Wang}, \bibinfo{author}{Y.~Gu}, \bibinfo{author}{S.~Zhang},
  \bibinfo{author}{L.-Q. Chen}, \bibinfo{title}{Flexoelectricity in solids:
  Progress, challenges, and perspectives}, \bibinfo{journal}{Progress in
  Materials Science} \bibinfo{volume}{106} (\bibinfo{year}{2019})
  \bibinfo{pages}{100570}.

\bibitem[{Arias et~al.(2022)Arias, Catalan, and Sharma}]{arias2022emancipation}
\bibinfo{author}{I.~Arias}, \bibinfo{author}{G.~Catalan},
  \bibinfo{author}{P.~Sharma}, \bibinfo{title}{The emancipation of
  flexoelectricity}, \bibinfo{year}{2022}.

\bibitem[{Mocci et~al.(2021)Mocci, Barcel{\'o}-Mercader, Codony, and
  Arias}]{mocci2021geometrically}
\bibinfo{author}{A.~Mocci}, \bibinfo{author}{J.~Barcel{\'o}-Mercader},
  \bibinfo{author}{D.~Codony}, \bibinfo{author}{I.~Arias},
  \bibinfo{title}{Geometrically polarized architected dielectrics with apparent
  piezoelectricity}, \bibinfo{journal}{Journal of the Mechanics and Physics of
  Solids} \bibinfo{volume}{157} (\bibinfo{year}{2021}) \bibinfo{pages}{104643}.

\bibitem[{Zhang et~al.(2016)Zhang, Zhang, Shen, and Chen}]{zhang2016enhancing}
\bibinfo{author}{C.~Zhang}, \bibinfo{author}{L.~Zhang},
  \bibinfo{author}{X.~Shen}, \bibinfo{author}{W.~Chen},
  \bibinfo{title}{Enhancing magnetoelectric effect in multiferroic composite
  bilayers via flexoelectricity}, \bibinfo{journal}{Journal of Applied Physics}
  \bibinfo{volume}{119}~(\bibinfo{number}{13}) (\bibinfo{year}{2016})
  \bibinfo{pages}{134102}.

\bibitem[{Zhang et~al.(2021)Zhang, Yan, Wang, and Shao}]{zhang2021ultrahigh}
\bibinfo{author}{M.~Zhang}, \bibinfo{author}{D.~Yan},
  \bibinfo{author}{J.~Wang}, \bibinfo{author}{L.-H. Shao},
  \bibinfo{title}{Ultrahigh flexoelectric effect of 3D interconnected porous
  polymers: modelling and verification}, \bibinfo{journal}{Journal of the
  Mechanics and Physics of Solids} \bibinfo{volume}{151} (\bibinfo{year}{2021})
  \bibinfo{pages}{104396}.

\bibitem[{McBride et~al.(2020)McBride, Davydov, and
  Steinmann}]{mcbride2020modelling}
\bibinfo{author}{A.~McBride}, \bibinfo{author}{D.~Davydov},
  \bibinfo{author}{P.~Steinmann}, \bibinfo{title}{Modelling the flexoelectric
  effect in solids: a micromorphic approach}, \bibinfo{journal}{Computer
  Methods in Applied Mechanics and Engineering} \bibinfo{volume}{371}
  (\bibinfo{year}{2020}) \bibinfo{pages}{113320}.

\bibitem[{Guinovart-Sanju{\'a}n et~al.(2019)Guinovart-Sanju{\'a}n, Merodio,
  L{\'o}pez-Realpozo, Vajravelu, Rodr{\'\i}guez-Ramos, Guinovart-D{\'\i}az,
  Bravo-Castillero, and Sabina}]{guinovart2019asymptotic}
\bibinfo{author}{D.~Guinovart-Sanju{\'a}n}, \bibinfo{author}{J.~Merodio},
  \bibinfo{author}{J.~C. L{\'o}pez-Realpozo}, \bibinfo{author}{K.~Vajravelu},
  \bibinfo{author}{R.~Rodr{\'\i}guez-Ramos},
  \bibinfo{author}{R.~Guinovart-D{\'\i}az},
  \bibinfo{author}{J.~Bravo-Castillero}, \bibinfo{author}{F.~J. Sabina},
  \bibinfo{title}{Asymptotic homogenization applied to flexoelectric rods},
  \bibinfo{journal}{Materials} \bibinfo{volume}{12}~(\bibinfo{number}{2})
  (\bibinfo{year}{2019}) \bibinfo{pages}{232}.

\bibitem[{Mohammadi et~al.(2014)Mohammadi, Liu, and
  Sharma}]{mohammadi2014theory}
\bibinfo{author}{P.~Mohammadi}, \bibinfo{author}{L.~Liu},
  \bibinfo{author}{P.~Sharma}, \bibinfo{title}{A theory of flexoelectric
  membranes and effective properties of heterogeneous membranes},
  \bibinfo{journal}{Journal of Applied Mechanics}
  \bibinfo{volume}{81}~(\bibinfo{number}{1}).

\bibitem[{Sharma et~al.(2007)Sharma, Maranganti, and
  Sharma}]{sharma2007possibility}
\bibinfo{author}{N.~Sharma}, \bibinfo{author}{R.~Maranganti},
  \bibinfo{author}{P.~Sharma}, \bibinfo{title}{On the possibility of
  piezoelectric nanocomposites without using piezoelectric materials},
  \bibinfo{journal}{Journal of the Mechanics and Physics of Solids}
  \bibinfo{volume}{55}~(\bibinfo{number}{11}) (\bibinfo{year}{2007})
  \bibinfo{pages}{2328--2350}.

\bibitem[{Geers and Yvonnet(2016)}]{geers2016multiscale}
\bibinfo{author}{M.~Geers}, \bibinfo{author}{J.~Yvonnet},
  \bibinfo{title}{Multiscale modeling of microstructure--property relations},
  \bibinfo{journal}{MRS Bulletin} \bibinfo{volume}{41}~(\bibinfo{number}{8})
  (\bibinfo{year}{2016}) \bibinfo{pages}{610--616}.

\bibitem[{Geers et~al.(2010)Geers, Kouznetsova, and
  Brekelmans}]{geers2010multi}
\bibinfo{author}{M.~G. Geers}, \bibinfo{author}{V.~G. Kouznetsova},
  \bibinfo{author}{W.~Brekelmans}, \bibinfo{title}{Multi-scale computational
  homogenization: Trends and challenges}, \bibinfo{journal}{Journal of
  computational and applied mathematics}
  \bibinfo{volume}{234}~(\bibinfo{number}{7}) (\bibinfo{year}{2010})
  \bibinfo{pages}{2175--2182}.

\bibitem[{Kaczmarczyk et~al.(2008)Kaczmarczyk, Pearce, and
  Bi{\'c}ani{\'c}}]{kaczmarczyk2008scale}
\bibinfo{author}{{\L}.~Kaczmarczyk}, \bibinfo{author}{C.~J. Pearce},
  \bibinfo{author}{N.~Bi{\'c}ani{\'c}}, \bibinfo{title}{Scale transition and
  enforcement of RVE boundary conditions in second-order computational
  homogenization}, \bibinfo{journal}{International Journal for Numerical
  Methods in Engineering} \bibinfo{volume}{74}~(\bibinfo{number}{3})
  (\bibinfo{year}{2008}) \bibinfo{pages}{506--522}.

\bibitem[{Kouznetsova et~al.(2004)Kouznetsova, Geers, and
  Brekelmans}]{kouznetsova2004multi}
\bibinfo{author}{V.~Kouznetsova}, \bibinfo{author}{M.~G. Geers},
  \bibinfo{author}{W.~Brekelmans}, \bibinfo{title}{Multi-scale second-order
  computational homogenization of multi-phase materials: a nested finite
  element solution strategy}, \bibinfo{journal}{Computer methods in applied
  Mechanics and Engineering} \bibinfo{volume}{193}~(\bibinfo{number}{48-51})
  (\bibinfo{year}{2004}) \bibinfo{pages}{5525--5550}.

\bibitem[{Lesi{\v{c}}ar et~al.(2014)Lesi{\v{c}}ar, Tonkovi{\'c}, and
  Sori{\'c}}]{lesivcar2014second}
\bibinfo{author}{T.~Lesi{\v{c}}ar}, \bibinfo{author}{Z.~Tonkovi{\'c}},
  \bibinfo{author}{J.~Sori{\'c}}, \bibinfo{title}{A second-order two-scale
  homogenization procedure using $C^{1}$ macrolevel discretization},
  \bibinfo{journal}{Computational mechanics}
  \bibinfo{volume}{54}~(\bibinfo{number}{2}) (\bibinfo{year}{2014})
  \bibinfo{pages}{425--441}.

\bibitem[{Lesi{\v{c}}ar et~al.(2016)Lesi{\v{c}}ar, Sori{\'c}, and
  Tonkovi{\'c}}]{lesivcar2016large}
\bibinfo{author}{T.~Lesi{\v{c}}ar}, \bibinfo{author}{J.~Sori{\'c}},
  \bibinfo{author}{Z.~Tonkovi{\'c}}, \bibinfo{title}{Large strain, two-scale
  computational approach using C1 continuity finite element employing a second
  gradient theory}, \bibinfo{journal}{Computer methods in applied mechanics and
  engineering} \bibinfo{volume}{298} (\bibinfo{year}{2016})
  \bibinfo{pages}{303--324}.

\bibitem[{Lesi{\v{c}}ar et~al.(2017)Lesi{\v{c}}ar, Tonkovi{\'c}, and
  Sori{\'c}}]{lesivcar2017two}
\bibinfo{author}{T.~Lesi{\v{c}}ar}, \bibinfo{author}{Z.~Tonkovi{\'c}},
  \bibinfo{author}{J.~Sori{\'c}}, \bibinfo{title}{Two-scale computational
  approach using strain gradient theory at microlevel},
  \bibinfo{journal}{International Journal of Mechanical Sciences}
  \bibinfo{volume}{126} (\bibinfo{year}{2017}) \bibinfo{pages}{67--78}.

\bibitem[{Lesi{\v{c}}ar(2015)}]{lesivcar2015multiscale}
\bibinfo{author}{T.~Lesi{\v{c}}ar}, \bibinfo{title}{Multiscale modeling of
  heterogeneous materials using second-order homogenization}, Ph.D. thesis,
  \bibinfo{school}{University of Zagreb. Faculty of Mechanical Engineering and
  Naval Architecture.}, \bibinfo{year}{2015}.

\bibitem[{Ghasemi et~al.(2017)Ghasemi, Park, and Rabczuk}]{ghasemi2017level}
\bibinfo{author}{H.~Ghasemi}, \bibinfo{author}{H.~S. Park},
  \bibinfo{author}{T.~Rabczuk}, \bibinfo{title}{A level-set based IGA
  formulation for topology optimization of flexoelectric materials},
  \bibinfo{journal}{Computer Methods in Applied Mechanics and Engineering}
  \bibinfo{volume}{313} (\bibinfo{year}{2017}) \bibinfo{pages}{239--258}.

\bibitem[{Nanthakumar et~al.(2017)Nanthakumar, Zhuang, Park, and
  Rabczuk}]{nanthakumar2017topology}
\bibinfo{author}{S.~Nanthakumar}, \bibinfo{author}{X.~Zhuang},
  \bibinfo{author}{H.~S. Park}, \bibinfo{author}{T.~Rabczuk},
  \bibinfo{title}{Topology optimization of flexoelectric structures},
  \bibinfo{journal}{Journal of the Mechanics and Physics of Solids}
  \bibinfo{volume}{105} (\bibinfo{year}{2017}) \bibinfo{pages}{217--234}.

\bibitem[{Deng et~al.(2017)Deng, Deng, Yu, and Shen}]{deng2017mixed}
\bibinfo{author}{F.~Deng}, \bibinfo{author}{Q.~Deng}, \bibinfo{author}{W.~Yu},
  \bibinfo{author}{S.~Shen}, \bibinfo{title}{Mixed finite elements for
  flexoelectric solids}, \bibinfo{journal}{Journal of Applied Mechanics}
  \bibinfo{volume}{84}~(\bibinfo{number}{8}).

\bibitem[{Mao et~al.(2016)Mao, Purohit, and Aravas}]{mao2016mixed}
\bibinfo{author}{S.~Mao}, \bibinfo{author}{P.~K. Purohit},
  \bibinfo{author}{N.~Aravas}, \bibinfo{title}{Mixed finite-element
  formulations in piezoelectricity and flexoelectricity},
  \bibinfo{journal}{Proceedings of the Royal Society A: Mathematical, Physical
  and Engineering Sciences} \bibinfo{volume}{472}~(\bibinfo{number}{2190})
  (\bibinfo{year}{2016}) \bibinfo{pages}{20150879}.

\bibitem[{Poya et~al.(2019)Poya, Gil, Ortigosa, and Palma}]{poya2019family}
\bibinfo{author}{R.~Poya}, \bibinfo{author}{A.~J. Gil},
  \bibinfo{author}{R.~Ortigosa}, \bibinfo{author}{R.~Palma}, \bibinfo{title}{On
  a family of numerical models for couple stress based flexoelectricity for
  continua and beams}, \bibinfo{journal}{Journal of the Mechanics and Physics
  of Solids} \bibinfo{volume}{125} (\bibinfo{year}{2019})
  \bibinfo{pages}{613--652}.

\bibitem[{He et~al.(2019)He, Javvaji, and Zhuang}]{he2019characterizing}
\bibinfo{author}{B.~He}, \bibinfo{author}{B.~Javvaji},
  \bibinfo{author}{X.~Zhuang}, \bibinfo{title}{Characterizing flexoelectricity
  in composite material using the element-free Galerkin method},
  \bibinfo{journal}{Energies} \bibinfo{volume}{12}~(\bibinfo{number}{2})
  (\bibinfo{year}{2019}) \bibinfo{pages}{271}.

\bibitem[{Yvonnet and Liu(2017)}]{yvonnet2017numerical}
\bibinfo{author}{J.~Yvonnet}, \bibinfo{author}{L.~Liu}, \bibinfo{title}{A
  numerical framework for modeling flexoelectricity and Maxwell stress in soft
  dielectrics at finite strains}, \bibinfo{journal}{Computer Methods in Applied
  Mechanics and Engineering} \bibinfo{volume}{313} (\bibinfo{year}{2017})
  \bibinfo{pages}{450--482}.

\bibitem[{Dasgupta and Sengupta(1990)}]{dasgupta1990higher}
\bibinfo{author}{S.~Dasgupta}, \bibinfo{author}{D.~Sengupta}, \bibinfo{title}{A
  higher-order triangular plate bending element revisited},
  \bibinfo{journal}{International Journal for Numerical Methods in Engineering}
  \bibinfo{volume}{30}~(\bibinfo{number}{3}) (\bibinfo{year}{1990})
  \bibinfo{pages}{419--430}.

\bibitem[{Maranganti et~al.(2006)Maranganti, Sharma, and
  Sharma}]{maranganti2006electromechanical}
\bibinfo{author}{R.~Maranganti}, \bibinfo{author}{N.~Sharma},
  \bibinfo{author}{P.~Sharma}, \bibinfo{title}{Electromechanical coupling in
  nonpiezoelectric materials due to nanoscale nonlocal size effects: Green’s
  function solutions and embedded inclusions}, \bibinfo{journal}{Physical
  Review B} \bibinfo{volume}{74}~(\bibinfo{number}{1}) (\bibinfo{year}{2006})
  \bibinfo{pages}{014110}.

\bibitem[{Sharma et~al.(2010)Sharma, Landis, and
  Sharma}]{sharma2010piezoelectric}
\bibinfo{author}{N.~Sharma}, \bibinfo{author}{C.~Landis},
  \bibinfo{author}{P.~Sharma}, \bibinfo{title}{Piezoelectric thin-film
  superlattices without using piezoelectric materials},
  \bibinfo{journal}{Journal of Applied Physics}
  \bibinfo{volume}{108}~(\bibinfo{number}{2}) (\bibinfo{year}{2010})
  \bibinfo{pages}{024304}.

\bibitem[{Abdollahi et~al.(2015)Abdollahi, Mill{\'a}n, Peco, Arroyo, and
  Arias}]{abdollahi2015revisiting}
\bibinfo{author}{A.~Abdollahi}, \bibinfo{author}{D.~Mill{\'a}n},
  \bibinfo{author}{C.~Peco}, \bibinfo{author}{M.~Arroyo},
  \bibinfo{author}{I.~Arias}, \bibinfo{title}{Revisiting pyramid compression to
  quantify flexoelectricity: A three-dimensional simulation study},
  \bibinfo{journal}{Physical Review B}
  \bibinfo{volume}{91}~(\bibinfo{number}{10}) (\bibinfo{year}{2015})
  \bibinfo{pages}{104103}.

\bibitem[{Shen and Hu(2010)}]{shen2010theory}
\bibinfo{author}{S.~Shen}, \bibinfo{author}{S.~Hu}, \bibinfo{title}{A theory of
  flexoelectricity with surface effect for elastic dielectrics},
  \bibinfo{journal}{Journal of the Mechanics and Physics of Solids}
  \bibinfo{volume}{58}~(\bibinfo{number}{5}) (\bibinfo{year}{2010})
  \bibinfo{pages}{665--677}.

\bibitem[{Hu and Shen(2010)}]{hu2010variational}
\bibinfo{author}{S.~Hu}, \bibinfo{author}{S.~Shen}, \bibinfo{title}{Variational
  principles and governing equations in nano-dielectrics with the flexoelectric
  effect}, \bibinfo{journal}{Science China Physics, Mechanics and Astronomy}
  \bibinfo{volume}{53}~(\bibinfo{number}{8}) (\bibinfo{year}{2010})
  \bibinfo{pages}{1497--1504}.

\bibitem[{Hill(1963)}]{hill1963elastic}
\bibinfo{author}{R.~Hill}, \bibinfo{title}{Elastic properties of reinforced
  solids: some theoretical principles}, \bibinfo{journal}{Journal of the
  Mechanics and Physics of Solids} \bibinfo{volume}{11}~(\bibinfo{number}{5})
  (\bibinfo{year}{1963}) \bibinfo{pages}{357--372}.

\bibitem[{Schr{\"o}der and Keip(2012)}]{schroder2012two}
\bibinfo{author}{J.~Schr{\"o}der}, \bibinfo{author}{M.-A. Keip},
  \bibinfo{title}{Two-scale homogenization of electromechanically coupled
  boundary value problems}, \bibinfo{journal}{Computational mechanics}
  \bibinfo{volume}{50}~(\bibinfo{number}{2}) (\bibinfo{year}{2012})
  \bibinfo{pages}{229--244}.

\bibitem[{Hrabok and Hrudey(1984)}]{hrabok1984review}
\bibinfo{author}{M.~Hrabok}, \bibinfo{author}{T.~Hrudey}, \bibinfo{title}{A
  review and catalogue of plate bending finite elements},
  \bibinfo{journal}{Computers \& Structures}
  \bibinfo{volume}{19}~(\bibinfo{number}{3}) (\bibinfo{year}{1984})
  \bibinfo{pages}{479--495}.

\bibitem[{Argyris et~al.(1968)Argyris, Fried, and Scharpf}]{argyris1968tuba}
\bibinfo{author}{J.~H. Argyris}, \bibinfo{author}{I.~Fried},
  \bibinfo{author}{D.~W. Scharpf}, \bibinfo{title}{The TUBA family of plate
  elements for the matrix displacement method}, \bibinfo{journal}{The
  Aeronautical Journal} \bibinfo{volume}{72}~(\bibinfo{number}{692})
  (\bibinfo{year}{1968}) \bibinfo{pages}{701--709}.

\bibitem[{Bell(1969)}]{bell1969refined}
\bibinfo{author}{K.~Bell}, \bibinfo{title}{A refined triangular plate bending
  finite element}, \bibinfo{journal}{International journal for numerical
  methods in engineering} \bibinfo{volume}{1}~(\bibinfo{number}{1})
  (\bibinfo{year}{1969}) \bibinfo{pages}{101--122}.

\bibitem[{Zervos et~al.(2001)Zervos, Papanastasiou, and
  Vardoulakis}]{zervos2001finite}
\bibinfo{author}{A.~Zervos}, \bibinfo{author}{P.~Papanastasiou},
  \bibinfo{author}{I.~Vardoulakis}, \bibinfo{title}{A finite element
  displacement formulation for gradient elastoplasticity},
  \bibinfo{journal}{International Journal for Numerical Methods in Engineering}
  \bibinfo{volume}{50}~(\bibinfo{number}{6}) (\bibinfo{year}{2001})
  \bibinfo{pages}{1369--1388}.

\bibitem[{Zervos et~al.(2009)Zervos, Papanicolopulos, and
  Vardoulakis}]{zervos2009two}
\bibinfo{author}{A.~Zervos}, \bibinfo{author}{S.-A. Papanicolopulos},
  \bibinfo{author}{I.~Vardoulakis}, \bibinfo{title}{Two finite-element
  discretizations for gradient elasticity}, \bibinfo{journal}{Journal of
  engineering mechanics} \bibinfo{volume}{135}~(\bibinfo{number}{3})
  (\bibinfo{year}{2009}) \bibinfo{pages}{203--213}.

\bibitem[{Manzari and Yonten(2013)}]{manzari2013c1}
\bibinfo{author}{M.~T. Manzari}, \bibinfo{author}{K.~Yonten},
  \bibinfo{title}{C1 finite element analysis in gradient enhanced continua},
  \bibinfo{journal}{Mathematical and Computer Modelling}
  \bibinfo{volume}{57}~(\bibinfo{number}{9-10}) (\bibinfo{year}{2013})
  \bibinfo{pages}{2519--2531}.

\bibitem[{Aifantis(2000)}]{aifantis2000gradient}
\bibinfo{author}{E.~C. Aifantis}, \bibinfo{title}{Gradient aspects of crystal
  plasticity at micro and macro scales}, in: \bibinfo{booktitle}{Key
  Engineering Materials}, vol. \bibinfo{volume}{177},
  \bibinfo{organization}{Trans Tech Publ}, \bibinfo{pages}{805},
  \bibinfo{year}{2000}.

\bibitem[{Krichen and Sharma(2016)}]{krichen2016flexoelectricity}
\bibinfo{author}{S.~Krichen}, \bibinfo{author}{P.~Sharma},
  \bibinfo{title}{Flexoelectricity: a perspective on an unusual
  electromechanical coupling}, \bibinfo{journal}{Journal of Applied Mechanics}
  \bibinfo{volume}{83}~(\bibinfo{number}{3}).

\bibitem[{Deng et~al.(2014)Deng, Liu, and Sharma}]{deng2014flexoelectricity}
\bibinfo{author}{Q.~Deng}, \bibinfo{author}{L.~Liu},
  \bibinfo{author}{P.~Sharma}, \bibinfo{title}{Flexoelectricity in soft
  materials and biological membranes}, \bibinfo{journal}{Journal of the
  Mechanics and Physics of Solids} \bibinfo{volume}{62} (\bibinfo{year}{2014})
  \bibinfo{pages}{209--227}.

\end{thebibliography}

\end{document}